\documentclass[10pt]{amsart}

\usepackage[utf8]{inputenc}
\usepackage[T1]{fontenc}

\usepackage{amsmath}
\usepackage{float}
\usepackage[all]{xy}
\usepackage{fancyhdr}
\usepackage{graphicx}

\usepackage{mathtools}
\usepackage{adjustbox}
\usepackage{enumitem}
\usepackage{stackrel}
\usepackage{scalerel}
\usepackage{euscript}
\usepackage{amsmath}
\usepackage{amsthm}
\usepackage{indentfirst}

\usepackage{mathpazo}
\usepackage{mathrsfs}
\usepackage{amsmath}
\usepackage{amssymb}
\usepackage{hyperref}
\usepackage{cleveref}
\usepackage{comment}
\usepackage{color}
\usepackage{tikz}
\usepackage{tikz-cd}
\usepackage{xparse}
\usetikzlibrary{matrix,arrows,decorations.pathmorphing}

\newtheorem{defn}{Definition}[section]
\newtheorem{thm}{Theorem}[section]
\newtheorem{prop}[thm]{Proposition}
\newtheorem{cor}[thm]{Corollary}
\newtheorem{lem}[thm]{Lemma}

\newtheorem{remark}[thm]{Remark}
\newtheorem{example}[thm]{Example}

\crefname{lem}{Lemma}{lemma}
\crefname{remark}{Remark}{remark}
\crefname{cor}{Corollary}{corollary}
\crefname{thm}{Theorem}{theorem}
\crefname{prop}{Proposition}{proposition}
\crefname{example}{Example}{example}
\crefname{defn}{Definition}{definition}
\crefname{notation}{Notation}{notation}
\crefname{appendix}{Appendix}{appendix}
\crefname{section}{Section}{section}

\newtheorem{innercustomthm}{Theorem}
\newenvironment{nnumthm}[1]
  {\renewcommand\theinnercustomthm{#1}\innercustomthm}
  {\endinnercustomthm}

\newcommand{\TTT}{\mathcal T}
\newcommand{\NNN}{\mathcal N}

\newcommand\Nb {\mathbb{N}}
\newcommand\Zb {\mathbb{Z}}

\newcommand\Cb {\mathbb{C}}

 \newcommand{\G}{\mathscr{G}}
\newcommand{\F}{\mathscr {F}}

\newcommand\Bs {\mathbf{s}}
\newcommand\Bt {\mathbf{t}}

\newcommand\CA {\EuScript{A}}
\newcommand\CB {\EuScript{B}}
\newcommand\CC {\EuScript{C}}

\newcommand\CH {\EuScript{H}}

\newcommand\CK {\EuScript{K}}
\newcommand\CL {\EuScript{L}}

\newcommand\CN {\EuScript{N}}

\newcommand\CS {\EuScript{S}}

\newcommand{\FA} {\mathfrak{A}}
\newcommand{\FB} {\mathfrak{B}}
\newcommand{\FC} {\mathfrak{C}}

\newcommand{\FH} {\mathfrak{H}}

\newcommand{\FX} {\mathfrak{X}}
\newcommand{\FY} {\mathfrak{Y}}

\newcommand{\Fo} {\mathfrak{o}}
\newcommand{\Ft} {\mathfrak{t}}

\usepackage{amsfonts} 
\DeclareMathSymbol{\shortminus}{\mathbin}{AMSa}{"39}
\usepackage{mathpazo}

\DeclareMathOperator{\ob}{Ob}
\DeclareMathOperator{\id}{Id}

\DeclareMathOperator{\cott}{\FC \Fo \Ft\Ft}

\newcommand{\mref}[1]{(\ref{#1})}

\DeclareMathOperator{\spn}{span}

\DeclareMathOperator{\Ran}{Ran}

\DeclareDocumentCommand{\lin}{m O{\cdot} O{\cdot}}{{}_{{#1}{ }} \langle #2, #3 \rangle}
\DeclareDocumentCommand{\rin}{m O{\cdot} O{\cdot}}{\langle #2, #3 \rangle_{{#1}{}} }
\DeclareDocumentCommand{\ot}{m}{\otimes_{{#1}{}} }
\DeclareDocumentCommand{\od}{m}{\odot_{{#1}{}} }
\DeclarePairedDelimiterX\braket[2]{\langle}{\rangle}{#1 \delimsize\vert #2}

\newcommand{\titleinfo}{Co-universal $C^{\ast}$-algebras for product systems over finite aligned subcategories of groupoids}
\newcommand{\titleinfoshort}{Co-universal $C^{\ast}$-algebras}
\newcommand{\authorinfo}{Feifei Miao, Liguang Wang, Wei Yuan}

\begin{document}

\title{\LARGE\textbf{\titleinfo}}
\author{\large\textsc{\authorinfo}}

\address{School of Mathematics and Information Science, Shandong Technology and Business University, Yantai, Shandong, 264005, China }
\email{mff100511@163.com }

\address{School of Mathematical Sciences, Qufu Normal University, Qufu, Shandong, 273165, China}
\email{wangliguang0510@163.com}

\address{AMSS, Chinese Academy of Sciences, Beijing, 100190,  China\\
and\\
School of Mathematical Sciences\\
University of Chinese Academy of Sciences, Beijing 100049, China}
\email{wyuan@math.ac.cn}

\begin{abstract}
The product systems over left cancellative small categories are introduced and studied in this paper. We also introduce the notion of compactly aligned product systems over finite aligned left cancellative small categories and its  Nica covariant  representations. The  existence of co-universal algebras for  injective, gauge-compatible, Nica covariant representations of compactly aligned product systems over finite aligned subcategories of groupoids is proved in this paper.

\end{abstract}

\subjclass[2010]{46L05}
\keywords{product system, left cancellative small category, Fell bundle, cosystem, $C^*$-envelope, co-universal $C^*$-algebra.}
 \thanks{Wang was supported in part by NSF of China (No. 11871303, No. 11971463, No. 11671133) and NSF of Shandong Province (No. ZR2019MA039, No. ZR2020MA008).  Yuan was supported in part by NSF of China (No. 11871303, No. 11871127, No. 11971463).}
\date{}
\maketitle

\section{Introduction}
In \cite{MR561974}, Cuntz and Krieger studied a class of C$^{\ast}$-algebras which are later called Cuntz-Krieger algebras. Since then, there have been many generalizations of Cuntz-Krieger algebras. Pimsner generalized Cuntz-Krieger algebras by considering C$^*$-algebras associated with C$^*$-correspondences \cite{MR1426840}. Motivated by Pimsner's work, Katsura further generalized Pimsner's construction and introduced a method to construct C$^{\ast}$-algebras from C$^{\ast}$-correspondences in \cite{MR2102572}. The algebras studied in \cite{MR2102572} are now called Cuntz-Pimsner algebras.

Every C$^{\ast}$-correspondence defines a product system over $\mathbb{N}$, i.e., a notion introduced by Fowler. In \cite{FowlerN2002}, Fowler initiated the study of C$^*$-algebras associated with product systems, and introduced a Cuntz-type algebra for a regular product system over a directed quasi-lattice ordered group. After that, for quite some time, the right notion of Cuntz algebras for product systems of non-injective $C^*$-correspondences is unclear. Eventually, Sims and Yeend made a breakthrough and gave the definition of Cuntz-Nica-Pimsner algebras for general product systems under some minor conditions \cite{Sims2010}. Later, Carlsen, Larsen, Sims, and Vittadello inroduced the co-universal Cuntz-Pimsner algebras of product systems. These algebras are co-universal for injective, Nica covariant, gauge-compatible representations of the product systems and can be viewed as a generalization of reduced crossed products of C$^*$-algebras. In \cite{MR2837016}, it is proved that the algebras satisfying such co-universal property do exist in many cases. In further work, Dor-On and Katsoulis showed that the $C^*$-envelope of the Fock tensor algebra $(\TTT^r_\FX)^+$ has the co-universal property for every compactly aligned product system $\FX$ over an abelian lattice ordered group (\cite{MR4053621}).  Motivated by this result, Dor-On, Kakariadis, Katsoulis, Laca, and Li characterized the co-universal algebras of compactly aligned product systems over group-embeddable right LCM semigroups by using C$^{\ast}$-envelopes for cosystems in \cite{DKKLL21}. Most recently, Shehnem generalized the work of  Dor-On, Kakariadis, Katsoulis, Laca, and Li and proved that the $C^*$-envelopes of $(\TTT^r_\FX)^+$ are  co-universal for injective compatible representations of product systems over arbitrary group-embeddable submonoids in \cite{envearb}.

In \cite{MR3549542}, Albandik and Meyer described  Cuntz-Pimsner algebras of product systems over  monoids from another perspective. Let $\CA$ and $\CB$ be $C^*$-algebras, and let $\cott(\CA, \CB)$ be the groupoid consisting of correspondences from $\CA$ to $\CB$ (i.e, $\CA$-$\CB$ correspondences) and $C^*$-correspondence isomorphisms.
The composition is defined by $(\FY, \FX) \mapsto \FX\otimes_{\CB} \FY:  \cott(\CB, \CC) \times \cott(\CA, \CB) \to \cott(\CA, \CC)$. Then $\cott$ is a bicategory where it's objects, arrows and 2 arrows are $C^*$-algebras, $C^*$-correspondences and $C^*$-correspondence isomorphisms, respectively. Albandik and Meyer described the  functor from a category to $\cott$ in detail and proved that if $P$ is a monoid, then the functor $P \to \cott$ can be identified with the product system over an opposite monoid of $P$. Moreover, they proved that
if $P \to \cott$ maps every object of $P$ to a regular correspondence, then
 the colimit of $P \to \cott$ can be regarded as the Cuntz-Pimsner algebra of the regular product system. Inspired by their results, we can introduce the concept of product systems over  left cancellative small categories.

A right LCM semigroup can be considered as a specific instance of left cancellative small categories, which has been extensively studied by B\'{e}dos, Kaliszewski, Quigg, and Spielberg in \cite{MR3909245} with regards to the Toeplitz and Cuntz-Krieger type algebras.
A natural question is whether there exist co-universal algebras related to the product systems over left cancellative small categories.


The purpose of this article is to introduce the concept of Nica covariant Toeplitz representations of compactly aligned product systems over finite aligned left cancellative small categories and
show that the technique developed in \cite{DKKLL21} can be used to characterize the co-universal algebras of product systems over finite aligned subcategories of groupoids.  More explicitly, we proved the following theorem.

\begin{nnumthm}{1}\label{thm:_main_result}
 Let $\FX = (\{\CA_{x}\}, \{\FX_{f}\})$ be a compactly aligned product system over a finite aligned subcategory of a groupoid $G$. Then the co-universal C$^*$-algebra for injective,  gauge-compatible, Nica covariant representations of $\FX$ exists. Moreover, the coaction of $G$ on the co-universal C$^*$-algebra is normal.
\end{nnumthm}
The paper is organized as follows.

In section 2, we introduce the concept of product systems over left cancellative small categories and study the representations on them. Furthermore, we introduce the concept of compactly aligned product systems over finite aligned left cancellative small categories and the concept of Nica covariant Toeplitz representations on these product systems. Moreover, we demonstrate \cref{Nica_cond}. In section 3, we extend the notion of cosystem introduced in \cite{DKKLL21} by replacing group with groupoid, and provide some basic results on cosystems. In section 4, we recall some basic facts of Fell bundles and graded $C^*$-algebras over discrete groupoids. In section 5, we give the proof of \cref{thm:_main_result}.

\section{Quasi-product systems over left cancellative small categories}
For a right Hilbert module $\FX$ over a C$^*$-algebra (we refer to \cite{L95} for the definition and basic properties of the right Hilbert $C^*$-modules), we use $\CL(\FX)$ to denote the C$^*$-algebra of all adjointable operators on $\FX$ and $\CK(\FX)$ to denote the C$^*$-algebra of all compact operators on $\FX$.

A $\CA$-$\CB$ $C^{\ast}$-correspondence of $C^{\ast}$-algebras $\CA$ and $\CB$ is a right Hilbert $\CB$-module $\FX$ equipped with a non-degenerate $*$-homomorphism $\phi : \CA \to \CL(\FX)$, i.e. $\{ \phi(a)\xi: a \in \CA, \xi \in \FX\}$ is dense in $\FX$. In the following, we will use $a \xi$ to denote the left action $\phi(a) \xi$.

\begin{defn}
A bimodule map $F:\FX \to \FY$ between two $\CA$-$\CB$ C$^*$-correspondences is called a C$^*$-correspondence map if $\rin{\CB}[\xi][\beta] = \rin{\CB}[F(\xi)][F(\beta)]$ for all $\xi$, $\beta \in \FX$. A $C^*$-correspondence map $F$ is called a \textit{C$^*$-correspondence isomorphism} if $F$ is a unitary.
\end{defn}

Let $\CA$, $\CB$ and $\CC$ be $C^{\ast}$-algebras. 
The interior tensor product $\FX \ot{\CB} \FY$ of a $\CA$-$\CB$ C$^{\ast}$-correspondence $\FX$ and a $\CB$-$\CC$ $C^*$-correspondence $\FY$ is the completion of the algebraic $\CB$-balanced tensor product $\FX \od{\CB} \FY$ with respect to the norm induced by the $\CC$-valued inner product defined by
\begin{align*}
	\rin{\CC}[\xi_1 \odot \eta_1][\xi_2 \odot \eta_2]:=  \rin{\CC}[\eta_1][\rin{\CB}[\xi_1][\xi_2]\eta_2], \quad \xi_1, \xi_2 \in \FX, \eta_1, \eta_2 \in \FY.
\end{align*}
It is well-known that $\FX \ot{\CB} \FY$ is a $\CA$-$\CC$ C$^*$-correspondence (see \cite{L95}).\\

In this paper, we only consider small categories. For a category $G$, we use $G(x,y)$ to denote the set of morphisms from $x$ to $y$, and $G(-, -)$ to denote the set of morphisms in $G$. The identity morphism of $x \in \ob(G)$ is written as $1_x$. For $f \in G(x,y)$, we use $\Bs(f)$ (resp. $\Bt(f)$) to denote the source (resp. target) of $f$, i.e., $\Bs(f) = x$ (resp. $\Bt(f) = y$). Let
\begin{align*}
	G(-,y) :&= \{f \in G(-, -): \Bt(f) = y\}, \quad  G(x, -) := \{f \in G(-, -): \Bs(f) = x\}.
\end{align*}
Recall that a category is left cancellative if $fg=fh$ implies $g = h$ for every (composable) morphisms in the category. For the rest of this section, we use $P$ to denote a left cancellative small category.

\subsection{Product systems over left cancellative categories and their Toeplitz representations}

\begin{defn}\label{def:_product_sys}
A quasi-product system over $P$ is a triple
   \begin{align*}
       \left (\{\CA_{x}\}_{x \in Ob(P)}, \{\FX_{f}\}_{f \in P(-,-)}, \{m_{g, f}\}_{\{(g,f) \in P(-,-)^2: \Bt(f) = \Bs(g)\}} \right ),
   \end{align*}
where $\CA_x$ is a  C$^*$-algebra, $\FX_{f}$ is a $\CA_{\Bt(f)}$-$\CA_{\Bs(f)}$ C$^*$-correspondence, and $m_{g, f}$ is a C$^*$-correspondence map from $\FX_g \ot{\CA_{\Bt(f)}} \FX_f$ to $\FX_{gf}$, subject to the following axioms
\begin{enumerate}
    \item $\FX_{1_x} = \CA_x$ for every $x \in \ob(P)$;
    \item $m_{1_x, f}(a \otimes_{\CA_x} \xi) = a \xi$ and $m_{f, 1_w}(\xi \otimes_{\CA_w} b) = \xi b$ for $f \in P(w,x)$, $a \in \CA_x$, $b \in \CA_w$, $\xi \in \FX_f$;
    \item for every $f \in P(w, x)$, $g \in P(x,y)$, $h \in P(y,z)$, the following diagram
    \begin{align*}
        \xymatrix @R=0.2in @C=0.6in{
            \FX_h \ot{\CA_y} \FX_g\ot{\CA_{x}} \FX_f \ar[d]^{I_{\FX_h} \ot{\CA_y} m_{g,f}}
            \ar[r]^-{m_{h,g} \ot{\CA_{x}} I_{\FX_f}} & \FX_{hg} \ot{\CA_{x}} \FX_f \ar[d]^{m_{hg,f}}\\
            \FX_h \ot{\CA_y}  \FX_{gf} \ar[r]^{m_{h,gf}} & \FX_{hgf}
        }
    \end{align*}
commutes.
\end{enumerate}
\end{defn}

\begin{remark}
A  quasi-product system over $P$ is amount to a strictly unitary lax functor from the opposite bicategory of $P$ to the bicategory of C$^*$-algebras, bimodules of C$^*$-algebras, and C$^*$-correspondence maps (see  \cite[Proposition 6.2]{MR3549542} and \cite{JY20} for more details).
\end{remark}

\begin{defn}
A quasi-product system
       \begin{align*}
       \left (\{\CA_{x}\}_{x \in Ob(P)}, \{\FX_{f}\}_{f \in P(-,-)}, \{m_{g, f}\}_{\{(g,f) \in P(-,-)^2: \Bt(f) = \Bs(g)\}} \right )
       \end{align*}
       over $P$ is called product system if every $m_{g, f}$ is  a C$^*$-correspondence isomorphism.
\end{defn}

In the following, we usually suppress $\{m_{g,f}\}$ from the notation of the product system, and abbreviate a product system over $P$ to $(\{\CA_{x}\}, \{\FX_{f}\})$, and write $m_{g, f}(\xi_g \otimes_{\CA_{\Bt(f)}}\eta_f)$ as $\xi_g \eta_f$.

From now on, we use $\FX = (\{\CA_{x}\}, \{\FX_{f}\})$ to denote a product system over $P$.

\begin{defn}\label{def:iso_rep_comp_alig_prod_system}
	A Toeplitz \textit{representation} of $\FX = (\{\CA_{x}\}, \{\FX_{f}\})$ into a $C^*$-algebra $\CB$, denoted by $\psi : \FX\to \CB$, is a family of linear maps $\{\psi_f: \FX_f \to \CB\}_{f \in P(-, -)}$ such that
	\begin{enumerate}
		\item $\psi_{1_x} : \CA_x \to \CB$ is a $*$-homomorphism for every $x \in \ob(P)$.
		\item $\psi_g(\beta) \psi_f(\xi) = \psi_{gf}(\beta \xi)$ for every $\beta \in \FX_g$, $\xi \in \FX_f$ and composable morphisms $g$, $f \in P(-, -)$.
		\item $\psi_f(\eta)^*\psi_f(\xi) = \psi_{1_{\Bs(f)}} (\rin{\CA_{\Bs(f)}}[\eta][\xi])$ for  $\xi, \eta \in \FX_f$, $f \in P(-,-)$.
	\end{enumerate}
	We say $\psi$ is  non-degenerate if the C$^*$-subalgebra $C^*(\psi)$ generated by $\{\psi_f(\FX_f)\}_{f \in P(-,-)}$ is non-degenerate, i.e., $\overline{\spn}\{ c \xi: c \in C^*(\psi), \xi \in \CB\}=  \CB$.
\end{defn}

\begin{remark}
	Let $\psi: \FX \to \CB$ be a Toeplitz representation of $\FX$. If $\rho$ is a $*$-homomorphism from $C^*(\psi)$ into a C$^*$-algebra $\CC$, then $\rho \circ \psi := \{ \rho \circ \psi_f: \FX_f \to \CC\}_{f \in P(-, -)}$ is also a Toeplitz representation of $\FX$.
\end{remark}

Let $\psi: \FX \to \CB$ be a Toeplitz representation of $\FX$. It is clear that $\psi_f$ is a contraction for every $ f \in P(-, -)$. Also $\psi_f$ is isometric if $\psi_{1_{\Bs(f)}}$ is injective. The Toeplitz representation $\psi$ is called injective if $\psi_{1_x}$ is injective for all $x \in \ob(P)$.

\begin{remark}
    Let $\FX=(\{\CA_{x}\}, \{\FX_{f}\})$ be a product system over $P$. For a $C^*$-algebra algebra $\CB$, let  $\Delta_{\CB}$ be the constant bifunctor $\Delta_{\CB}: P \to \cott$. More explicitely, $\Delta_{\CB}$ sends every object/morphism to $\CB$ (see Proposition 4.1.21 in \cite{JY20}).
    A strong transformation from  $\FX$ to $\Delta_{\CB}$ such that every component 1-morphism in $\cott(\CA_x, \CB)$ regarded as a right Hilbert $\CB$-module equals $\CB$  is a Toeplitz \textit{representation} of $\FX$ in $\CB$ (see \cite{MR3549542}).
\end{remark}

\begin{example} \label{Sub:_fock}
	For a quasi-product system $\FX = (\{\CA_{x}\}, \{\FX_{f}\})$, let
	\begin{align*}
		C_c(\FX):=\{\xi =(\xi_f)_{f \in P(-, -)}: \xi_f \in \FX_f, \mbox{the support of $(\xi_f)$ is a finite set}\}.
	\end{align*}
	Then $C_c(\FX)$ is a pre-Hilbert right module over the C$^*$-algebra $C_0(\{\CA_x\})$, where $C_0(\{\CA_x\})$ is the completion of
	\begin{align*}
		\{a:=(a_x)_{x \in \ob(P)}: a_x \in \CA_x, \mbox{and the support of $(a_x)$ is a finite set} \}
	\end{align*}
	with respect to the supremum norm. More explicitely, the right $C_0(\{\CA_x\})$-action on $C_c(\FX)$ and the $C_0(\{\CA_x\})$-valued inner product are given by
	\begin{align*}
		\xi \cdot a := \left (\xi_f \cdot a_{\Bs(f)} \right),\quad
		\rin{C_0(\CA_x)}[\xi][\eta] := \left (\sum_{f \in P(x, -)} \rin{\CA_{x}}[\xi_f][\eta_f] \right).
	\end{align*}
	Let $C_0(\FX)$ be the completion of $C_c(\FX)$.
	
	For every $\eta_g  \in \FX_g$, let  $L_{\eta_g}: C_0(\FX) \to C_0(\FX)$ be the operator defined by
	\begin{align*}
		(L_{\eta_g}\xi)_f =
		\begin{cases}
			\eta_g \xi_h, &\mbox{there exists $h$ such that $f = gh$};\\
			0, & \mbox{otherwise}.
		\end{cases}
	\end{align*}
Now we assume that $\FX$ is a product system. We see that the operator $L_{\eta_g}$ is adjointable. Indeed, the adjoint $L^*_{\eta_g}$ is given by
	\begin{align*}
		(L^*_{\eta_g}\xi)_f :=
		\begin{cases}
			m_{1_{\Bs(g)}, f} \left (\rin{\CA_{\Bs(g)}}[\eta_g][-] \otimes_{\CA_{\Bs(g)}} I \right )  m^*_{g, f} \xi_{gf}, & \text{ if $\Bt(f) = \Bs(g)$};\\
			0, & \mbox{otherwise}.
		\end{cases}
	\end{align*}
It is easy to check that $\{\eta_g \mapsto L_{\eta_g}: \FX_g \to \CL(C_0(\FX))\}_{g \in P(-,-)}$ is an injective Toeplitz representation of the product system $\FX$ (Note that $\|L_{\eta_g}\| = \|\eta_g\|$). This representation is called the \textit{Fock representation} of $\FX$. The norm closed algebra generated by  $\{L_{\eta_f}: \eta_f \in \FX_f, f \in P(-, -)\}$  is called the \textit{Fock tensor algebra} $(\TTT^r_\FX)^+$. And the \textit{Fock algebra} $\TTT^r_\FX$ is the $C^*$-algebra generated by $(\TTT^r_\FX)^+$.
\end{example}

\subsection{Compactly aligned product systems over finite aligned left cancellative categories}
Let $f \in P(-,-)$. We use $fP$ to denote the set $\{ f l : l \in P(-, -), \Bs(f)= \Bt(l)\}$.
\begin{defn}[\cite{MR3909245}]\label{def:_finitesubset}
	If there is a finite subset $\F$ of $P(-, -)$ such that $fP \cap gP = \cup_{\omega \in \F} \omega P$ for every $f, g \in P(-, -)$, then we say $P$ is   finite aligned.
\end{defn}

Note that the finite set $\F$ in \cref{def:_finitesubset} can be an empty set.
Recall that a quasi-lattice ordered group is a pair $(G, P)$, where $G$ is a discrete group with unit $e$ and $P$ is a subsemigroup of $G$ such that $P \cap P^{-1} = \{e\}$ and every finite subset of $G$ with an upper bound in $P$ has the least upper bound in $P$. It is obvious that $P$ is a finite aligned subcategory of $G$.
\begin{example}
	Let $(G', P')$ be a quasi-lattice ordered group.
	\begin{enumerate}
		\item Let $G$ be the category with $\ob(G) = \Zb^{+}$ and $G(i, j) = G'$ for every $i, j \in \ob(G)$. And $P$ is the wide subcategory of $G$, i.e., $\ob(P) = \ob(G)$, such that
		\begin{align*}
			P({i},{j}) =
			\begin{cases}
				P', &  \text{if $i=1,j>1$ or $i = j$;} \\
				\emptyset, & \text{otherwise.}
			\end{cases}
		\end{align*}
		Then $P$ is a finite aligned subcategory of $G$.
		
		\item The principle groupoid $G$ of $(G', P')$ is the category such that $\ob(G) = \{ x\}_{x\in P'}$ and $G(x, y)$ contains a unique morphism for every $x, y \in P'$. Let $P$ be the wide subcategory of $G$ with
		\begin{align*}
			P(x, y)=
			\begin{cases}
				G(x, y), &  \text{if } \quad y \le x;\\
				\emptyset, &  \text{ otherwise} .
			\end{cases}
		\end{align*}
		Then $P$ is a  finite aligned subcategory of $G$.
		
		\item Let $X$ be a $G'$-set, $G'$ acts on $X$. The transformation groupoid $G$ is the category with $\ob(G) = X$ and $G(x, y) = \{g\in G': g \cdot x =y\}$ (see example 8.1.15 in \cite{Sims20101}). Let $P$ be the wide subcategory of $G$ such that $P(x, y) = \{ g \in P': g \cdot x =y\}$. Then $P$ is a  finite aligned subcategory of $G$.
	\end{enumerate}
\end{example}

From now on, we assume that $P$ is finite aligned unless explicitly stated otherwise.

Let $\{ f_1, \cdots, f_n\} \subseteq P(-, -)$. It is not hard to check that there exists a finite subset of $P(-, -)$, denoted by $\vee\{ f_1, \cdots, f_n\}$ or $f_1 \vee  \cdots \vee f_n$, such that
\begin{enumerate}
	\item $f_1P \cap  f_2P \cdots \cap f_nP =\cup_{\omega \in \vee\{ f_1, \cdots, f_n\}} \omega P$,
	\item $h \notin k P$ for every two non-equal elements $h, k \in \vee \{ f_1, \cdots, f_n\}$.
\end{enumerate}
Moreover, if $\F$ is another finite subset of $P(-, -)$ satisfying the conditions $(1)$ and $(2)$, then every $h \in \F$ can be written as $h  = \omega k$, where $\omega \in \vee \{ f_1, \cdots, f_n\}$ and $k$ is an invertible morphism in $P(-, -)$.

Let $\FX = (\{\CA_{x}\}, \{\FX_{f}\})$ be a product system over $P$.
For $f, g \in P(-, -)$, we use $l^g_f$ to denote the $*$-homomorphism from $\CL(\FX_f)$ to $\CL(\FX_g)$ defined as follows
\begin{align*}
	l^{g}_f(A) :=\begin{cases}
		m_{f, h} (A \otimes_{\CA_{\Bs(f)}} I) m^*_{f, h}, & \mbox{if}\quad g = fh;\\
		0, & \text{otherwise},
	\end{cases} \quad A \in \CL(\FX_f).
\end{align*}
Note that $\CK(\FX_{1_{\Bt(g)}}) = \CA_{\Bt(g)}$, and $ l^g_{1_{\Bt(g)}}(a)(\xi_g) = a \xi_g$ for $\xi_g \in \FX_g$, $a \in \CA_{\Bt(g)}$ (see condition (2) in \cref{def:_product_sys}).

In \cite[Definition 5.7]{FowlerN2002}, Fowler introduced the notion of compactly aligned product
systems over quasi-lattice ordered groups. This notion is extended to product systems over right LCM semigroups by Brownlowe, Larsen, and Stammeier in \cite{MR3900375}. We now extend this notion to product systems over finite aligned left cancellative categories.

\begin{defn}
A product system $\FX = (\{\CA_{x}\}, \{\FX_{f}\})$ over $P$ is called \textit{compactly aligned} if $ l^{\omega}_f(S) l^{\omega}_g(T) \in \CK(\FX_{\omega})$ for every $S \in \CK(\FX_f)$, $T \in \CK(\FX_g)$, $f, g \in P(-, -)$, $\omega \in f \vee g$.
\end{defn}

For the rest of this section, we assume that $\FX = (\{\CA_{x}\}, \{\FX_{f}\})$ is a compactly aligned product system over $P$.
Let $\psi : = \{ \psi_f : \FX_f \to \CB\}_{f \in P(-, -)}$ be a Toeplitz representation of $\FX$. For every $f \in P(-, -)$, there is a $*$-homomorphism $\psi^{(f)} : \CK(\FX_f) \to \CB$ such that
	\begin{align}\label{eq:compact_decomposable}
		\psi^{(f)}(\theta_{\xi, \eta}) = \psi_{f}(\xi)\psi_{f}(\eta)^*,
	\end{align}
	where $\theta_{\xi, \eta}(-) := \xi\rin{\CA_{\Bs(f)}}[\eta][-]$ is the rank one operator defined by $\xi, \eta \in \FX_f$ (see \cite[p. 202]{MR1426840},  \cite[Lemma 2.2]{MR1658088} and \cite[Remark 1.7]{MR1722197}).
Let  $\CB(\FH)$ be the algebra of bounded operators on a Hilbert space $\FH$.  If $\CB \subseteq \CB(\FH)$, then it is not hard to check that
\begin{align*}
	\overline{\psi^{(f)}(\CK(\FX_f)) \FH} =  \overline{\spn}\{ \psi_f(\xi_f) \zeta: \xi_f \in \FX_f, \zeta \in \FH\}.
\end{align*}
By Proposition 2.5 in \cite{L95}, there is a strict-SOT (here SOT denotes strong operator topology) continuous extension $\CL(\FX_f) \to \CB(\FH)$, which is also denoted by $\psi^{(f)}$. For  an approximate unit $ \{ e_{i}\}_{i \in \mathbb{I}}$ of $\CK(\FX_f)$, $\psi^{(f)}(I) = \lim_{i} \psi^{(f)}(e_i) (strict-SOT)$ is the projection from $\FH$ onto $\overline{\psi^{(f)}(\CK(\FX_f)) \FH}$.

\begin{lem}\label{lem:_Topelitz_app}
Let $\psi: \FX \to \CC$ be a Toeplitz representation. For $g \in P(-, -)$, let $ \{ e_{i}\}_{i \in \mathbb{I}} $ be an increasing approximate unit of $\CK(\FX_g)$.
	Then
	\begin{align}\label{lem:appcom}
		\lim_i \psi^{(g)}(e_i) \psi^{(\omega)}(S) = \psi^{(\omega)}(S),
	\end{align}
	where $S \in \CK(\FX_\omega)$, $\omega = g k$, $k \in P(-, -)$.
\end{lem}

\begin{proof}
	Since $\FX_g \FX_k = \FX_\omega$ and rank one compact operators are linear dense in $\CK(\FX_\omega)$, we only need to show that
	\begin{align*}
		\lim_i \psi^{(g)}(e_i) \psi^{(\omega)}(\theta_{\xi_g\eta_k, \zeta_\omega}) = \psi^{(\omega)}(\theta_{\xi_g\eta_k, \zeta_\omega}), \quad  \forall \xi_g \in \FX_g, \eta_k \in \FX_k, \zeta_\omega \in \FX_\omega.
	\end{align*}
	Note that $\lim_i e_i \xi_g = \xi_g$, we have
	\begin{align*}
		\lim_i \psi^{(g)}(e_i)\psi^{(\omega)}(\theta_{\xi_g\eta_k, \zeta_\omega}) &=\lim_i \psi^{(g)}(e_i)\psi_g(\xi_g) \psi_k(\eta_k) \psi_{\omega}(\zeta_\omega)^*\\
		&=\lim_i \psi_g(e_i \xi_g) \psi_k(\eta_k) \psi_{\omega}(\zeta_\omega)^* =\psi^{(\omega)}(\theta_{\xi_g\eta_k, \zeta_\omega}).
	\end{align*}
\end{proof}

\begin{lem}\label{lem_the_reg_lift}
	Let $\psi : \FX \to \CB(\FH)$ be a  Toeplitz representation. Then
	\begin{align*}
		\psi^{(g)}(T) \psi^{(\omega)} (I) =  \psi^{(\omega)}\Big(l_g^{\omega}(T)\Big),
	\end{align*}
	for $T \in \CL(\FX_g)$, $\omega = g v$, $g, v \in P(-, -)$.
\end{lem}

\begin{proof}
	We see that
	\begin{align*}
		\psi^{(\omega)}\Big(l_g^{\omega}(T)\Big) \psi_g(\eta_g) \psi_{v}(\zeta_{v}) \xi &= \psi_\omega \Big(l_g^{\omega}(T) \eta_g \zeta_{v}\Big) \xi \\
		&= \psi_g(T \eta_g) \psi_{v}(\zeta_{v}) \xi = \psi^{(g)}(T) \psi_g(\eta_g) \psi_{v}(\zeta_{v}) \xi,
	\end{align*}
	for $\eta_g \in \FX_g$, $\zeta_{v} \in \FX_{v}$, $\xi \in \FH$.  Since $\spn\{ \psi_g(\eta_g) \psi_{v} (\zeta_{v}) \xi: \eta_g \in \FX_g, \zeta_{v} \in \FX_{v}, \xi \in \FH\}$ is dense in $\psi_{\omega}(\FX_{\omega}) \FH$, the lemma is proved.
\end{proof}

\subsection{The Nica covariant representation of compactly aligned product systems}

\begin{defn}
	Let $\psi$ be a Toeplitz representation of $\FX$ and $\rho: C^*(\psi) \to \CB(\CH)$ is a $*$-representation.  The Toeplitz representation $\psi$ is called \textit{Nica covariant with respect to $\rho$} if
	\begin{align*}
		( \rho \circ \psi)^{(f)}(I)\ ( \rho \circ \psi)^{(g)}(I) = \begin{cases}
			\vee_{h \in f \vee g} ( \rho \circ \psi)^{(h)}(I), & \text{ if   $f \vee g \neq \emptyset$}; \\
			0, & \text{otherwise},
		\end{cases}
	\end{align*}
	for every $f, g \in P(-, -)$. If $\psi$ is Nica covariant with respect to every $*$-representation, then $\psi$ is called \textit{Nica covariant}.
\end{defn}

\begin{remark}
	If $\psi$ is Nica covariant, then $\rho \circ \psi$ is Nica covariant for every $*$-homomorphism $\rho: C^*(\psi) \to \CC$.
\end{remark}

\begin{remark}
	Assume that $\psi$ is Nica covariant with respect to a $*$-representation $\rho: C^*(\psi) \to \CB(\FH)$. Then we have
	\begin{align*}
		(  \rho\circ \psi)^{(f)}(I)( \rho\circ  \psi)^{(g)}(I) = (\rho\circ  \psi)^{(g)}(I)(\rho \circ \psi)^{(f)}(I), \quad \forall f, g \in P(-, -).
	\end{align*}
\end{remark}

\begin{prop}\label{lem:_Nicacovariant_vertical}
	Let $\psi $ be a Toeplitz representation of $\FX$. Assume that $\psi$ is Nica covariant with respect to an  injective $*$-representation $\rho$. Let $f, g \in P(-, -)$. If $\Bs(f) \neq \Bt(g)$, then
	\begin{align*}
		\psi_f(\xi) \psi_g(\eta) = 0, \quad  \forall \xi \in \FX_f, \eta \in \FX_g.
	\end{align*}
\end{prop}

\begin{proof}
	Since $1_{\Bs(f)} \vee 1_{\Bt(g)} = \emptyset$, we have $(\rho\circ \psi)^{(1_{\Bs(f)})}(I) (\rho\circ \psi)^{(1_{\Bt(g)})}(I) = 0$. Thus $\psi_{1_{\Bs(f)}}(a) \psi_{1_{\Bt(g)}}(b) = 0$ for every $a \in \CA_{\Bs(f)}, b\in \CA_{\Bt(g)}$. Note that $\FX_f \cdot \CA_{\Bs(f)}$ (resp. $\CA_{\Bt(g)} \cdot \FX_g$) is dense in $\FX_f$ (resp. $\FX_g$). We have $\psi_f(\xi) \psi_g(\eta) = 0$, for every $\xi \in \FX_f, \eta \in \FX_g$.
\end{proof}
\begin{lem}\label{lem:_incluP}
	Let $\psi: \FX \to \CB(\FH)$ be a Toeplitz representation. Assume that $\psi$ is Nica covariant with respect to $\id_{C^*(\psi)}$. For every $f, g \in P(-, -)$, $\xi_f \in \FX_f$, we have
	\begin{align}
		\psi_f(\xi_f)\psi^{(g)} (I) &=
		\begin{cases}
			\psi^{(fg)}(I)\psi_{f}(\xi_f), & \Bs(f) = \Bt(g);\\
			0, & \Bs(f) \neq \Bt(g),
		\end{cases} \label{equ:lem:_incluP_1}\\
		\psi^{(g)}(I) \psi_f(\xi_f) &=
		\begin{cases}
			\psi_f(\xi_f) \Big(\vee_{v_i} \psi^{(v_i)}(I) \Big), &  g \vee f = \{ f v_i\}^n_{i= 1};\\
			0, &  g\vee f = \emptyset.
		\end{cases}\label{equ:lem:_incluP_2}
	\end{align}
	In particular, $\psi^{(f)}(S_f) \psi^{(g)}(I) =  \psi^{(g)}(I) \psi^{(f)}(S_f)$ for every $S_f \in \CK(\FX_f)$.
\end{lem}

\begin{proof}
	By \cref{lem:_Nicacovariant_vertical}, $\psi_f(\xi_f) \psi^{(g)} (I) = 0$ if $\Bs(f) \neq \Bt(g)$. Assume that $\Bs(f) = \Bt(g)$. We have
	\begin{align*}
		\psi_f(\xi_f) \psi^{(g)} (I)\left(\psi_g(\eta_{g}) \zeta \right) & =\psi_{f g}(\xi_f \eta_{g}) \zeta = \psi^{(fg)}(I)\psi_{f}(\xi_f)\Big(\psi_g( \eta_{g}) \zeta \Big),
	\end{align*}
	for every $\eta_g \in \FX_g$, $\zeta \in \FH$.
	
	If $g\vee f = \emptyset$, then $\psi^{(g)}(I) \psi^{(f)}(I) =0$. Thus $\psi^{(g)}(I) \psi_f(\xi_f) = 0$. Assume that $g \vee f = \{ f v_i\}^n_{i= 1}$. We have
	\begin{align*}
		\psi^{(g)}(I) \psi_f(\xi_f) &= \psi^{(g)}(I) \psi^{(f)}(I) \psi_f(\xi_f) \\
		&= \Big(\vee_{f v_i \in g \vee f} \psi^{(f v_i)}(I) \Big) \psi_f (\xi_f) = \psi_f(\xi_f) \Big(\vee_{v_i} \psi^{(v_i)}(I) \Big).
	\end{align*}
	It follows from \cref{equ:lem:_incluP_1} and \cref{equ:lem:_incluP_2}   that
	\begin{align*}
		\psi^{(f)}( \theta_{\xi_f, \eta_f }) \psi^{(g)}(I) =  \psi^{(g)}(I) \psi^{(f)}( \theta_{\xi_f, \eta_f }), \quad \forall \xi_f, \eta_f \in \FX_f.
	\end{align*}
	Then $\psi^{(f)}(S_f) \psi^{(g)}(I) =  \psi^{(g)}(I) \psi^{(f)}(S_f)$ for every $S_f \in \CK(\FX_f)$.
\end{proof}

\begin{remark}\label{remarkpro}
	With the notations and assumption as in \cref{lem:_incluP}. For every $f \in P(-, -)$ and a finite subset $\F$ of $P(-, -)$, there is a finite subset $\G \in P(-, -)$ such that
	\begin{align*}
		\psi_f(\xi_f) \Big( \Pi_{g \in \F} \psi^{(g)}(I) \Big ) = \Big ( \Pi_{h \in \G} \psi^{(h)}(I) \Big ) \psi_f(\xi_f), \quad \forall \xi_f \in \FX_f.
	\end{align*}
	And $\Pi_{g \in \F} \psi^{(g)}(I) \psi_f(\xi_f)$ is in the span of
	\begin{align*}
		\{\psi_f(\xi_f) \Pi_{h \in \G} \psi^{(h)}(I): \mbox{$\G$ is a finite subset of $P(-, -)$}\}.
	\end{align*}
\end{remark}

\begin{lem}\label{lem:sot_lem}
	With the notations and assumption as in \cref{lem:_incluP}. For every $T_f \in \CK(\FX_f)$, $T_g \in \CK(\FX_g)$, we have
	\begin{align*}
		\lim_{i} \psi^{(g)}(e_i) \psi^{(f)}(T_f) \psi^{(g)}(T_g) = \psi^{(f)}(T_f) \psi^{(g)}(T_g),
	\end{align*}
	where $ \{ e_{i}\}_{i \in \mathbb{I}}$ is an approximate unit of $\CK(\FX_g)$.
\end{lem}

\begin{proof}
	Assume that $f\vee g= \{\omega_1, \cdots, \omega_n\}$. By \cref{lem_the_reg_lift} and \cref{lem:_incluP}, we have
	\begin{equation}\label{eq:left_eqq}
		\begin{aligned}
			&\psi^{(f)}(T_f)\psi^{(g)}(T_g)\\
			=& \psi^{(f)}(T_f)\psi^{(g)}(T_g) \psi^{(f)}(I) \psi^{(g)}(I) \\
			=&  \psi^{(f)}(T_f)\psi^{(g)}(T_g) \Big( \psi^{(\omega_1)}(I) +\sum_{ 2\leq j \leq n} \psi^{(\omega_j)}(I)\big(I- \vee_{k=1}^{j-1}\psi^{(\omega_k)}(I)\big) \Big )\\
			=& \sum_{i= 1}^n \psi^{(\omega_i)} \left (l_f^{\omega_i}(T_f) l_g^{\omega_i}(T_g) \right) - \sum_{ 2\leq j \leq n} \psi^{(\omega_j)}\Big(l_f^{\omega_j}(T_f) l_g^{\omega_j}(T_g)\Big) \Big (\vee_{k=1}^{j-1} \psi^{(\omega_k)}(I) \Big).
		\end{aligned}
	\end{equation}
	Then \cref{lem:_Topelitz_app} implies the result.
\end{proof}

\begin{thm}\label{Nica_cond}
	Let $\psi$ be a Toeplitz representation of a compactly aligned product system $\FX$ over a finite aligned left cancellative small category. Then $\psi$ is Nica covariant if and only if there is an injective $*$-representation $\rho: C^*(\psi) \to \CB(\FH)$ such that $\psi$ is Nica covariant with respect to $\rho$.
\end{thm}

\begin{proof}
	We only need to show that $\psi$ is Nica covariant with respect to every $*$-homomorphism $\rho: C^*(\psi) \to \CB(\FH_\rho)$ if the Toeplitz representation $\psi: \FX \to \CB(\FH)$ is Nica covariant with respect to $\id_{C^*(\psi)}$.
	
	By \cref{lem:sot_lem}, we have
	\begin{align*}
		(\rho \circ \psi)^{(g)}(I) \rho(\psi^{(f)}(T_f)) \rho(\psi^{(g)}(T_g)) = \rho(\psi^{(f)}(T_f)) \rho(\psi^{(g)}(T_g)),
	\end{align*}
	for every $f, g \in P(-, -)$, $T_f \in \CK(\FX_f)$, and $T_g \in \CK(\FX_g)$. This implies that
	\begin{align*}
		(\rho \circ \psi)^{(f)}(T_f) (\rho \circ \psi)^{(g)}(I) =  (\rho \circ \psi)^{(g)}(I) (\rho \circ \psi)^{(f)}(T_f).
	\end{align*}
	Thus $(\rho \circ \psi)^{(f)}(I) (\rho \circ \psi)^{(g)}(I) = (\rho \circ \psi)^{(g)}(I) (\rho \circ \psi)^{(f)}(I)$. Moreover, $(\rho \circ \psi)^{(f)}(T_f) (\rho \circ \psi)^{(g)}(I) = 0$ if $f \vee g = \emptyset$.
	
	Assume that $f \vee g= \{ \omega_1, \cdots, \omega_n\}$. Since
	\begin{align*}
		\psi^{(f)}(I) \psi^{(g)}(I) \leq \sum_{j} \psi^{(\omega_j)}(I),
	\end{align*}
	we have
	\begin{align*}
		\psi^{(g)}(T_g^*)  \psi^{(f)}(T_f^*) \psi^{(f)}(T_f)\psi^{(g)}(T_g) \leq \sum_{j}\psi^{(\omega_j)} \Big(l^{\omega_j}_g(T^*_g) l^{\omega_j}_f(T^*_f) l^{\omega_j}_f(T_f) l^{\omega_j}_g(T_g)\Big).
	\end{align*}
	This implies that
	\begin{align*}
		(\rho \circ \psi)^{(f)}(T_f)(\rho \circ \psi)^{(g)}(T_g)  = (\rho \circ \psi)^{(f)}(T_f)(\rho \circ \psi)^{(g)}(T_g)\Big(\vee_{i=1}^n(\rho \circ \psi)^{(\omega_i)}(I)\Big).
	\end{align*}
	Since $\vee_{i=1}^{n}(\rho \circ \psi)^{(\omega_i)}(I)$ is a sub-projection of $(\rho \circ \psi)^{(f)}(I)(\rho \circ \psi)^{(g)}(I)$, we have
	\begin{align*}
		(\rho \circ \psi)^{(f)}(I)(\rho \circ \psi)^{(g)}(I) =  \vee_{i=1}^{n}(\rho \circ \psi)^{(\omega_i)}(I).
	\end{align*}
\end{proof}

\begin{example}\label{lem:Fockr}
	With the notation as in \cref{Sub:_fock},  and  let $\pi : C_0(\{ \CA_x\}_{x \in \ob(P)}) \to \CB(\FH_\pi)$ be a non-degenerate injective $*$-representation. Then $C_0(\FX) \ot{C_0(\{ \CA_x\})} \FH_\pi$ is a Hilbert space with the inner product defined by
	\begin{align*}
		\rin{\Cb}[\xi \otimes \zeta][ \xi' \otimes \zeta'] = \rin{\Cb}[\zeta][\pi(\rin{C_0(\{ \CA_x\})}[\xi][\xi']) \zeta'], \quad \xi, \xi' \in C_0(\FX), \zeta, \zeta' \in \FH_\pi.
	\end{align*}
	Since $\pi$ is injective, the $*$-homomorphism $A \mapsto A \otimes I : \CL(C_0(\FX)) \to \CB(C_0(\FX) \ot{C_0(\{ \CA_x\})} \FH_\pi)$ is injective (see Chapter 4 in \cite{L95}). In particular, $\eta_f \mapsto \widetilde{L}_{\eta_f} := L_{\eta_f} \otimes I: \FX_f \to \CB(C_0(\FX) \ot{C_0(\{ \CA_x\})} \FH_\pi)$ is an injective Toeplitz representation. Note that
	\begin{align*}
		\Ran( \widetilde{L}^{(f)}(I)) = \overline{\spn} \{ \xi_h \otimes \zeta: h \in fP, \xi_h \in \FX_h, \zeta \in \FH_\pi\},
	\end{align*}
	for every $f \in P(-, -)$. Then, it is not hard to see that $\widetilde{L}$ is Nica covariant (with respect to $\id_{C^*(\widetilde{L})}$). By \cref{Nica_cond}, the Fock representation $L$ is Nica covariant.
\end{example}

\begin{thm}\label{thm:_Nica-Topelitz}
	There exists an injective Nica covariant Toeplitz representation $i_\FX = \{ i_{\FX, f} : \FX_f \to \CB(\FH)\}_{f \in P(-, -)}$ of $\FX$ such that for every Nica covariant Toeplitz representation $\psi$ of $\FX$, there exists a (unique) $\ast$-homomorphism $\psi': C^*(i_\FX) \to C^*(\psi)$ rendering the following diagram
	\begin{align*}
		\xymatrix @R=0.2in {
			&  \FX \ar[ld]_{i_\FX} \ar[rd]^-{\psi} &  \\
			C^*(i_\FX) \ar[rr]^-{\psi'} & &C^*(\psi)
		}
	\end{align*}
	commutative.
\end{thm}

\begin{proof}
	A Nica covariant Toeplitz representation $\psi : \FX \to \CB(\FH)$ is said to be cyclic if the $C^*$-algebra $C^*(\psi)$ admits a cyclic vector in $\FH$. It is not hard to check that every Nica covariant Toeplitz representation $\psi : \FX \to \CB(\FH)$ can be decomposed as a direct sum of cyclic Nica covariant Toeplitz representations.
	
	Let $\CS$ be the set of cyclic Nica covariant Toeplitz representation of $\FX$. By \cref{lem:Fockr}, we have
	\begin{align*}
		i_{\FX, f} = \oplus_{\psi \in \CS} \psi_f: \FX_f \to \CB(\oplus_{\psi \in \CS} \FH_{\psi})
	\end{align*}
	is an injective Nica covariant Toeplitz representation satisfying the universal property.
\end{proof}

In the following, we use $\NNN\TTT_\FX$ to denote the C$^*$-algebra $C^*(i_\FX)$. By the universal property of $i_\FX$, $\NNN\TTT_\FX$ is the unique up to isomorphism, and is called the \textit{Nica-Toeplitz algebra} of $\FX$. And the \textit{Nica-Toeplitz tensor algebra} $\NNN\TTT_\FX^+$ of $\FX$ is the norm closed subalgebra of $\NNN\TTT_\FX$ generated by
\begin{align*}
	\{ i_{\FX,f}(\xi) : \xi \in \FX_f, f \in P(-, -)\}.\newline
\end{align*}

The following lemma will be used in the proof of \cref{thm:_main_result}.
\begin{lem} \label{lem:the_op}
	Let $\psi: \FX \to \CB(\FH)$ be a Nica covariant Toeplitz representation of $\FX$, and $f, g \in P(-, -)$. If $f \vee g =  \emptyset$, then $\psi_f(\xi_f)^* \psi_g(\eta_g) = 0$ for every $\xi_f \in \FX_f$, $\eta_g \in \FX_g$. Assume that $f \vee g = \{\omega_1 = fu_1 = gv_1, \ldots, \omega_n = f u_n = g v_n\}$. Then
	\begin{align*}
		\psi_f(\xi_f)^* \psi_g(\eta_g) \in \overline{\spn}\{ &\psi_{u_i}(\xi_{u_i}) \psi_{v_i}(\eta_{v_i})^* \Pi_{k \in \F} \psi^{(k)}(I): \\
		&\xi_{u_i} \in \FX_{u_i}, \eta_{v_i} \in \FX_{v_i}, \text{$\F$ is a finite subset of $P(-, -)$,  $i = 1, \ldots, n$}  \}.
	\end{align*}
\end{lem}

\begin{proof}
	If $f \vee g =  \emptyset$, then $\psi_f(\xi_f)^* \psi_g(\eta_g)  = \psi_f(\xi_f)^*\psi^{(f)}(I) \psi^{(g)}(I) \psi_g(\eta_g) = 0$.
	
	Assume that $f \vee g = \{\omega_1= fu_1 = gv_1, \ldots, \omega_n = f u_n = g v_n\}$. For $S \in \CK(\FX_f)$, $T \in \CK(\FX_g)$, we have
	\begin{equation}
		\begin{aligned}\label{eq:_one}
			\psi_f(\xi_f)^* \psi^{(f)}(S) \psi^{(g)}(T)\psi_g(\eta_g)  =& \psi_f(\xi_f)^* \psi^{(f)}(S) \psi^{(f)}(I) \psi^{(g)}(I) \psi^{(g)}(T) \psi_g(\eta_g)  \\
			=& \psi_f(\xi_f)^* \psi^{(f)}(S) \Big[ \vee_{\omega_i \in f\vee g}\psi^{(\omega_i)}(I) \Big] \psi^{(g)}(T)\psi_g(\eta_g) \\
			=&  \psi_f(\xi_f)^* \psi^{(f)}(S) \Big[\sum_{i} \psi^{(\omega_i)}(I) - \sum_{1 \leq i < j \leq n} \psi^{(\omega_i)}(I) \psi^{(\omega_j)}(I) \\
			&\quad + \sum_{1 \leq i < j < k \leq n} \psi^{(\omega_i)}(I) \psi^{(\omega_j)}(I) \psi^{(\omega_k)}(I) + \cdots  \\
			& \quad+ (-1)^{n-1} \psi^{(\omega_1)}(I) \cdots \psi^{(\omega_{n})}(I)\Big]\psi^{(g)}(T) \psi_g(\eta_g).
		\end{aligned}
	\end{equation}
	By \cref{lem_the_reg_lift} and the fact that  $\FX$ is compactly aligned, we have
	\begin{align*}
		\psi_f(\xi_f)^* \psi^{(f)}(S) & \psi^{(\omega_i)}(I) \psi^{(g)}(T)\psi_g(\eta_g) = \psi_f(\xi_f)^*\psi^{(\omega_i)}(l^{\omega_i}_f(S) l^{\omega_i}_g(T))\psi_g(\eta_g)   \\
		&  \in  \overline{\spn}\{  \psi_f(\xi_f)^* \psi_{\omega_i}(\zeta_{\omega_i}) \psi_{\omega_i}(\zeta'_{\omega_i})^* \psi_g(\eta_g) :  \zeta_{\omega_i}, \zeta'_{\omega_i} \in \FX_{\omega_i} \}.
	\end{align*}
	Since $\FX_{\omega_i} = \FX_f \FX_{u_i} = \FX_g \FX_{v_i}$, by the condition (3) of \cref{def:iso_rep_comp_alig_prod_system}, we have
	\begin{align*}
		\psi_f(\xi_f)^* \psi^{(f)}(S) \psi^{(\omega_i)}(I) \psi^{(g)}(T)\psi_g(\eta_g) \in \overline{\spn}\{\psi_{u_i}(\xi_{u_i}) \psi_{v_i}(\eta_{v_i})^*: \xi_{u_i} \in \FX_{u_i}, \eta_{v_i} \in \FX_{v_i}\}.
	\end{align*}
	Recall that (see \cref{lem:_incluP})
	\begin{align}\label{eq:_two}
		\psi^{(\omega_i)}(I) \psi^{(g)}(T) \psi_g(\eta_g)= \psi^{(g)}(T)\psi_g(\eta_g) \psi^{(v_i)}(I).
	\end{align}
	Then \cref{eq:_one} and \cref{eq:_two} implies that $\psi_f(\xi_f)^* \psi^{(f)}(S) \psi^{(g)}(T) \psi_g(\eta_g)$ is in
	\begin{align*}
		\overline{\spn}\{ \psi_{u_i}(\xi_{u_i}) \psi_{v_i}(\eta_{v_i})^* \Pi_{k \in \F} \psi^{(k)}(I):  &\xi_{u_i} \in \FX_{u_i}, \eta_{v_i} \in \FX_{v_i}, \\
		& \text{ $\F$ is a finite subset of $P(-, -)$}, i = 1, \ldots, n \}.
	\end{align*}
	Since $\CK(\FX_h)\FX_h$ is dense in $\FX_h$ for all $h \in P(-, -)$, the lemma is proved.
\end{proof}

For a Nica covariant Toeplitz representation $\psi : \FX \to \CB(\FH)$ of $\FX$, let
\begin{align}\label{gene:algebra}
	[C^*(\psi)] :=\overline{\spn}\{ &\psi_f(\xi) \psi_g(\eta)^* \Pi_{k \in \F} \psi^{(k)}(I): \nonumber \\
	& f, g \in P(-, -), \xi \in \FX_f, \eta \in \FX_g,  \text{$\F$ is a finite subset of $P(-, -)$ }  \}.
\end{align}
Since
\begin{align*}
	C^*(\psi) = \overline{\spn}\{&\psi_{f_1}(\xi_{f_1}) \psi_{g_1}(\eta_{g_1})^* \cdots \psi_{f_n}(\xi_{f_n}) \psi_{g_n}(\eta_{g_n})^* : \\
	& \xi_{f_i} \in \FX_{f_i}, \eta_{g_i} \in \FX_{g_i}, f_i, g_i \in P(-,-), 1 \leq i \leq n \in \Nb^{+}\},
\end{align*}
we have $C^*(\psi) \subset [C^*(\psi)]$ by \cref{lem:_Nicacovariant_vertical}, \cref{lem:_incluP}, and \cref{lem:the_op}.

\section{Co-systems and their $C^{\ast}$-envelopes}

\begin{defn}
 A $*$-representation of a left cancellative small category $P$ is a map $\pi$ from $P(-,-)$ to a $C^*$-algebra such that
\begin{enumerate}
    \item $\pi(f)\pi(g) = \delta_{\Bs(f), \Bt(g)} \pi(fg)$;
    \item $\pi(f)^* \pi(f) = \pi(1_{\Bs(f)})$;
    \item  $\pi(f) \pi(f)^* = \pi(1_{\Bt(f)})$
\end{enumerate}
for every $f, g \in P(-, -)$.
\end{defn}
A discrete groupoid is a category in which every morphism is invertible. For the rest of this section, we use $G$ to denote a discrete groupoid.

\begin{remark} Let $\CB$ be a $C^*$-algebra.
We see that $\pi : G \to \CB$ is a $*$-representation if and only if
\begin{enumerate}
    \item $\pi(f)\pi(g) = \delta_{\Bs(f), \Bt(g)} \pi(fg)$;
    \item $\pi(f)^* = \pi(f^{-1})$
\end{enumerate}
for every $f, g \in G(-, -)$.
\end{remark}

	Let $l^2(G)$ be the Hilbert space with the orthonormal basis $\{e_k\}_{k \in G(-, -)}$. The left regular $*$-representation $\lambda : G \to \CB(l^2(G))$ is the $*$-representation defined by
	\begin{align*}
		\lambda_g(e_k) := \begin{cases}
			e_{gk}, & \text{if}\quad \Bs(g) = \Bt(k); \\
			0, & \text{otherwise}.
		\end{cases}
	\end{align*}
	The \textit{reduced groupoid $C^{\ast}$-algebra $C^{\ast}_r(G)$ of $G$} is the C$^*$-algebra generated by the image of $\lambda$ (see \cite{Sims20101}). And the \textit{full groupoid $C^*$-algebra $C^*(G)$} is the universal C$^*$-algebra induced by all the $*$-representations of $G$. In particular, there is a $*$-homomorphism, which is also denoted by $\lambda$, from $C^*(G)$ to $C^*_r(G)$ defined by $U_g \mapsto \lambda_g$, where $U_g$ is the generator of $C^*(G)$ associated with $g$ (see Chapter 9 in \cite{Sims20101}).

\begin{remark}\label{rem:subalgebra}
Let $\pi: G \to \CB(\FH)$ be a $*$-representation of $G$. For every $x \in \ob(G)$, $\pi(1_x)$ is a projection and $f \mapsto \pi(f)|_{\pi(1_x) \FH} : G(x, x) \to \CB(\pi(1_x) \FH)$ is a unitary representation of the group $G(x, x)$. Conversely, by an argument similar to the one used in the proof of Stinespring's dilation theorem, it can be shown that every unitary representation $\pi_0: G(x,x) \to \CB(\CK)$ can be promoted to a $*$-representation of $G$. Indeed, we define a Hermitian form on the vector space $\FH_0 :=  \spn \{f \otimes \xi: \xi \in \CK, f \in G(x, -)\}$ as follows:
    \begin{align}\label{equ:subalgebra_1}
        \braket{f \otimes \xi}{g \otimes \beta}:= \delta_{\Bt(f), \Bt(g)} \braket{\xi}{\pi_0(f^{-1}g)\beta}.
    \end{align}
It is routine to check that the Hermitian form defined above is a semi-inner product. Then the Hermitian form \mref{equ:subalgebra_1} defines an inner product on the quotient space $\FH_0 / \CN$, where $\CN:=\{\zeta \in \FH_0: \braket{\zeta}{\zeta} = 0\}$. Let $\FH$ be completion of inner product space $\FH_0/ \CN$, and $\pi$ be the $*$-representation defined as follows:
    \begin{align*}
        \pi(h)(f \otimes \xi) := \delta_{\Bs(h), \Bt(f)} (hf) \otimes \xi, \quad \forall h \in G(-, -).
    \end{align*}
    It is clear that $\CK$ can be canonically identified with the subspace $\{1_x \otimes \xi : \xi \in \CK\} \subset \FH$. More explicitly, let $W$ be the isometry defined by $\xi \mapsto 1_x \otimes \xi: \CK \to \FH$. Then $W^*(f \otimes \xi) = \delta_{\Bt(f), x} \pi_0(f)\xi$, and $W\pi_0(h)W^* = \pi(h)$ for every $h \in G(x,x)$.
\end{remark}

Recall that an operator algebra is a norm closed subalgebra of $\CB(\FH)$ (see \cite{MR2111973, MR1976867}).
Let $\FA_1 \subset \CB(\FH_1)$, $\FA_2 \subset \CB(\FH_2)$ are two operator algebras. We use $\FA_1 \otimes \FA_2$ to denote the minimal tensor product of $\FA_1$ and $\FA_2$. Recall that $\FA_1 \otimes \FA_2$ is the completion of the algebraic tensor product $\FA_1\odot \FA_2$ of $\FA_1$ and $\FA_2$ with respect to the norm induced by the operator norm on $\CB(\FH_1 \otimes \FH_2)$ (see \cite{MR1065833}).

\begin{defn}
	A \textit{coaction} of $G$ on an operator algebra $\FA$ is a completely isometric homomorphism $\delta : \FA \to \FA \otimes C^*(G)$ such that $\spn \{\cup_{g \in G(-, -)} \FA^\delta_g\}$ is norm-dense in $\FA$, where
	\begin{align*}
		\FA^\delta_g := \{ a \in \FA:  \delta(a) = a \otimes U_g \}.
	\end{align*}
	If the map $(\id \otimes \lambda) \delta$ is injective, then the coaction $\delta$ is called \textit{normal}.
\end{defn}

\begin{example}\label{expl:coaction_full_G_alg}
By the universal property of $ C^*(G)$, there is a $*$-homomorphism
    \begin{align*}
        \delta_G : C^*(G) \to C^*(G) \otimes C^*(G)
    \end{align*}
such that $\delta_G(U_g) = U_g \otimes U_g$. By Remark \ref{rem:subalgebra}, we can identify $U_{1_x} C^*(G) U_{1_x}$ with $C^*(G(x,x))$ for $x \in \ob(G)$. Note that the restriction of $\delta_G$ to $C^*(G(x, x))$ is the canonical coaction $\delta_{G(x, x)}$ of $C^*(G(x, x))$ (see example 4.2.2 in \cite{MR2397671}). Since $\delta_{G(x, x)}$ is injective (see \cite{DKKLL21}), we have $\delta_G$ is injective. In particular, $\delta_G$ is a coaction of $G$ on $C^*(G)$. It is easy to check that $\delta_G$ is normal if and only if $\delta_{G(x, x)}$ is normal for every $x \in \ob(G)$.
\end{example}

\begin{remark}\label{rem:_algebra}
    Let $\delta: \FA \to \FA \otimes C^*(G)$ be an coaction. Then $\delta$ satisfies the coaction identity
\begin{align*}
(\delta \otimes \id_{C^*(G)})\circ \delta = (\id_{\FA} \otimes \delta_G)\circ \delta.
\end{align*}
If $\FA$ is a $C^*$-algebra, then $\delta$ is an injective $*$-homomorphism by Proposition 2.11 and Exercise 2.1 in \cite{MR1976867}. Moreover,
\begin{align*}
(\FA^\delta_g)^* = \{ a^* \in \FA:  \delta(a^*) = a^* \otimes U_{g^{-1}}\} = \FA^\delta_{g^{-1}}.
\end{align*}
\end{remark}

\begin{defn}
	A \textit{cosystem} is a pair $(\FA, \delta)$, where $\FA$ is an operator algebra, and $\delta : \FA \to \FA\otimes C^*(G)$ is a coaction. A morphism $(\FA, \delta) \to (\FA', \delta')$ between cosystems is an equivariant completely contractive homomorphism $\phi : \FA \to \FA'$, i.e., $\delta' \phi = (\phi \otimes \id) \delta$.
\end{defn}

The following proposition is a generalization of Fell absorption principle.
\begin{prop}[Generalized Fell absorption principle]\label{prop:_Generalized Fell absorption principle}
Let $\FX= (\{ \CA_x\}, \{ \FX_f\})$ be a quasi-product system over a left cancellative small category $P$.  Let $\pi$ be a $*$-representation of $P$ in $\CL(\FY)$, where $\FY$ is a right Hilbert module. Then there exists a partial isometry $W : C_0(\FX) \otimes \FY \to C_0(\FX) \otimes \FY $ such that $L_{\eta_g} \otimes \pi(g)  = W(L_{\eta_g}  \otimes I) W^*$, where $C_0(\FX) \otimes \FY$ is the exterior tensor product of $C_0(\FX)$ and $\FY$, $g \in P(-, -)$.
\end{prop}
\begin{proof}
It is obvious that we can regard $\FX_f$ as a subspace of  $C_0(\FX)$.        Let $W : C_0(\FX) \otimes \FY \to C_0(\FX) \otimes \FY $ be the operator defined as follows
\begin{align*}
W (\eta_f \otimes \xi ):= \eta_f \otimes \pi(f)\xi,
\end{align*}
    where $ \eta_f \in \FX_f$, $\xi \in \FY$. It is not hard to check that $W^*( \eta_f \otimes \xi )= \eta_f \otimes \pi(f)^*\xi$,
\begin{align*}
W^*W(\eta_f \otimes \xi) = \eta_f \otimes \pi(1_{\Bs(f)})\xi, \quad W W^*(\eta_f \otimes \xi) = \eta_f \otimes \pi(1_{\Bt(f)}) \xi,
\end{align*}
and $W (L_{\eta_g} \otimes I)W^* (\eta_f \otimes \xi) = (L_{\eta_g} \otimes \pi(g))(\eta_f \otimes \xi)$.
\end{proof}

The following fact is an immediate consequence of \cref{prop:_Generalized Fell absorption principle} and its proof is omitted (see, for example, proposition 3.4 in \cite{DKKLL21} for more details).

\begin{lem}\label{pro:_delta_tran_norm}
    Let $\delta_r : \FA \to \FA \otimes C^*_r(G)$ be a completely isometric homomorphism, where $\FA$ is an operator algebra. If $\spn\{\cup_{g \in G(-, -)} \{a \in \FA: \delta_r(a) = a \otimes \lambda_g\} \} $ is norm-dense in $\FA$, then $\FA$ admits a normal coaction $\delta$ of $G$ such that $\delta (b) = b \otimes U_g$ for $b \in \{ a \in \FA: \delta_r(a) = a \otimes \lambda_g\}$, $g \in G(-, -)$. In particular, $\delta_r = (\id \otimes \lambda) \delta$.
\end{lem}

\begin{remark}\label{exa:_reduced_groupoid}
By \cref{prop:_Generalized Fell absorption principle}, the reduced groupoid $C^*$-algebra $C^*_r(G)$ admits a $*$-homomorphism
\begin{align*}
\delta^r_G : C^*_r(G) \to C^*_r(G) \otimes C^*_r(G)
\end{align*}
    such that $\delta^r_G(\lambda_g) = \lambda_g \otimes \lambda_g$. It is not hard to see that $\delta^r_G$ is injective. Then there exists a normal coaction of $\delta: C^*_r(G) \to C^*_r(G) \otimes C^*(G)$ such that $\delta^r_G = (\id \otimes \lambda) \delta$ by \cref{pro:_delta_tran_norm}.
\end{remark}

A \textit{$C^*$-cover} of an operator algebra $\FA$ is a pair $(\CB, \iota)$, where $\CB$ is a C$^*$-algebra and $\iota : \FA \to \CB$ is a completely isometric homomorphism such that $\CB$ is generated by the image of $\iota$. The \textit{$C^*$-envelope} of $\FA$ is the terminal object $(C^*_{e}(\FA), \iota)$ in the category of C$^*$-covers of $\FA$. More explicitely, the terminality condition means that for every $C^*$-cover $(\CC, \iota')$ of $\FA$ there exists a (necessarily unique and surjective) $*$-homomorphism $\phi: \CC \to C^*_{e}(\FA)$ rendering the following diagram
\begin{align*}
	\xymatrix @R=0.2in {
		&  \FA \ar[ld]_{\iota'} \ar[rd]^-{\iota} &  \\
		\CC \ar[rr]_-{\phi}  & &C^*_{e}(\FA)
	}
\end{align*}
commutative (see \cite{MR253059}, \cite{MR566081},  \cite{ MR2132691}).

\begin{defn}
    A $C^*$-cover of a cosystem $(\FA, \delta)$ is a triple $(\CC, \delta_\CC ,\iota)$, where $\CC$ is a C$^*$-algebra, $\delta_\CC : \CC \to \CC \otimes C^*(G)$ is a coaction, and $\iota: \FA \to \CC$ is an equivariant completely contractive homomorphism such that $(\CC, \iota)$ is a C$^*$-cover of $\FA$.

The $C^{\ast}$-envelope of $(\FA, \delta)$ is a $C^{\ast}$-cover $(C^{\ast}_{e}(\FA, \delta), \delta_{e},\iota_{e})$ of $(\FA, \delta)$ satisfying the following property: for every $C^{\ast}$-cover $(\CC', \delta', \iota')$ of $(\FA, \delta)$, there exists an equivariant ${\ast}$-homomorphism $\phi : \CC' \to C^{\ast}_{e}(\FA, \delta)$ rendering the following diagram
  \begin{align*}
\xymatrix @R=0.2in {
    &  \FA \ar[ld]_{\iota'} \ar[rd]^-{\iota_{e}} &  \\
   \CC' \ar[rr]_-{\phi}  & &C^*_{e}(\FA, \delta)
   }
\end{align*}
commutative.
\end{defn}

In the following, we often omit $\delta_{e}$ and $\iota_{e}$ and write the C$^*$-envelope of $(\FA, \delta)$ simply as $C^*_{e}(\FA, \delta)$. The following fact is known to experts, at least when $G$ is a group (see Theorem 3.8 in \cite{DKKLL21}). For the sake of completeness, we sketch the proof.

\begin{prop}\label{lem:_envelope_for_the_cosystem}
	Let $(\FA, \delta)$ be a cosystem and $(C^*_{e}(\FA), \iota)$ be the C$^*$-envelope of $\FA$. Then
	\begin{align}\label{equ:envelope_for_the_cosystem_1}
		\big(C^*(\{\iota(a_g) \otimes U_g : a_g \in \FA^\delta_g,  g \in G(-, -)\}), \id \otimes \delta_G, (\iota \otimes \id_{C^*(G)}) \delta \big)
	\end{align}
	is the $C^*$-envelope of the cosystem $(\FA, \delta)$.
\end{prop}

\begin{proof}
    The homomorphism $(\iota \otimes \id_{C^*(G)}) \circ \delta : \FA \to  C^*_{e}(\FA) \otimes C^*(G)$ is completely isometric, since the minimal tensor product of completely isometric homomorphisms is completely isometric (see Theorem 2.2 in \cite{MR1065833}). It is also clear that $\id \otimes \delta_G$ is a coaction of $G$ on the C$^*$-algebra $C^*(\{\iota(a_g) \otimes U_g : a_g \in \FA^\delta_g,  g \in G(-, -)\})$ ($\subset C^*_e(\FA) \otimes C^*(G)$). Thus the triple \mref{equ:envelope_for_the_cosystem_1} is a C$^*$-cover of $(\FA, \delta)$.

For a $C^*$-cover $(\CC', \delta', \iota')$, there exists a $*$-homomorphism $\phi : \CC' \to C^*_{e}(\FA)$ such that $\iota = \phi \circ \iota'$. Let $\psi : = (\phi \otimes \id_{C^*(G)}) \delta'$.
It is obvious that $\psi$ is an equivariant  $*$-homomorphism such that $(\iota \otimes \id_{C^*(G)}) \delta = \psi \circ \iota'$.
\end{proof}

Let $\FA$ be an operator algebra and $\delta_r : \FA \to \FA \otimes C^*_r(G)$ be a completely isometric homomorphism satisfying the conditions of \cref{pro:_delta_tran_norm}. Then there exists a normal coaction $\delta: \FA \to \FA \otimes C^*(G)$ such that $\delta_r = (\id \otimes \lambda) \delta$. Let $\overline{\delta}_{C^*_r(G)} : C^*_r(G) \to C^*_r(G) \otimes C^*(G)$ be the normal coaction described in \cref{exa:_reduced_groupoid}. Then
\begin{align}\label{cor:_normal_coaction_relation}
    (C^*(\{ \iota(a_g) \otimes \lambda_g: a_g \in \FA^{\delta}_g, g \in G(-, -)\}), \id \otimes \overline{\delta}_{C^*_r(G)}, (\iota \otimes \id_{C^*_r(G)}) \delta_r)
\end{align}
is the $C^*$-envelope of the cosystem $(\FA, \delta)$, where $\iota$ is the completely isometric homomorphism from $\FA$ to $C^*_{e}(\FA)$. And the coaction $ \id \otimes \overline{\delta}_{C^*_r(G)}$ is normal (see the proof of Corollary 3.9 in \cite{DKKLL21} for more details).

\section{Fell bundles over groupoids and C$^*$-algebras graded by groupoids}

Fell bundles over groups were first introduced and studied by Fell \cite{MR0457620}, under the notion of $C^{\ast}$-algebraic bundle. It is a powerful tool in the study of graded $C^{\ast}$-algebras, and many well-known $C^*$-algebras are naturally graded. We briefly recall some basic facts of Fell bundles over discrete groupoids, and refer the readers to \cite{MR2881538} and \cite{Kumjian1998} for a detailed discussion of Fell bundles over groupoids.

Let $G$ be a discrete groupoid and $\FB = \{\FB_f\}_{f \in G(-, -)}$ be a Banach bundle over $G$ (see \cite[Definition 13.4 in Chapter II]{MR936628}). To simplify the notation, we use $\FB_x$ to denote $\FB_{1_x}$. A multiplication on $\FB$ is a family of continuous bilinear maps
\begin{align*}
  \{\FB_f \times \FB_g \to \FB_{fg}\}_{(f,g) \in G(-, -)^2, \Bs(f) = \Bt(g)}
\end{align*}
such that
\begin{enumerate}
    \item $(\xi_f \eta_g) \zeta_h = \xi_f (\eta_g \zeta_h)$ whenever the multiplication is defined,
    \item $\| \xi_f \eta_g\| \leq \|\xi_f\| \| \eta_g\|$ for all $\xi_f \in \FB_f$, $\eta_g \in \FB_g$, $f, g \in G(-, -)$, $\Bs(f) = \Bt(g)$.
\end{enumerate}
An involution on $\FB$ is a family of continuous involutive conjugate linear maps
\begin{align*}
    \xi_f \mapsto \xi_f^* : \FB_{f} \to \FB_{f^{-1}}, \quad f \in G(-, -).
\end{align*}

\begin{defn}[\cite{Kumjian1998}]\label{def:_the_fell_bundle_groupoid}
A Fell bundle over a discrete groupoid $G$ is a Banach bundle $\FB$ over $G$ equipped with a multiplication and an involution such that
\begin{enumerate}
    \item $(\xi_f \eta_g)^* = \eta^*_g \xi^*_f$ for $\xi_f \in \FB_f$, $\eta_g \in \FB_g$, $f, g \in G(-, -)$, $\Bs(f) = \Bt(g)$,
 \item $\| \xi_f^* \xi_f\| = \| \xi_f\|^2$ for $\xi_f \in \FB_f$, $f \in G(-, -)$,
    \item $\xi_f^* \xi_f$ is a positive element of the C$^*$-algebra $\FB_{\Bs(f)}$ for every $f \in G(-,-)$ and $\xi_f \in \FB_f$.
    \end{enumerate}
\end{defn}

For the rest of this subsection, we assume $\FB = \{\FB_f\}_{f \in G(-, -)}$ is a Fell bundle over $G$.
\begin{remark}
    For every $f \in G(-,-)$, $\FB_f$ together with the $\FB_{\Bs(f)}$-valued inner product $\rin{\FB_{\Bs(f)}}[\xi_f][\eta_f] := \xi^*_f \eta_f$ and the $\FB_{\Bt(f)}$-valued inner product $\lin{\FB_{\Bt(f)}}[\xi_f][\eta_f] = \xi_f \eta_f^*$ is a Hilbert $\FB_{\Bt(f)}$-$\FB_{\Bs(f)}$-bimodule, i.e., $\lin{\FB_{\Bt(f)}}[\xi_f][\eta_f] \zeta_f = \xi_f \rin{\FB_{\Bs(f)}}[\eta_f][\zeta_f]$. And $\FB^*_f = \FB_{f^{-1}}$.
\end{remark}

We see that $\FB$ is a quasi-product system, where the $C^*$-correspondence map $m_{f, g} : \FB_f \ot{\FB_{s(f)}} \FB_g \to \FB_{fg}$ is defined by the multiplication in $\FB$.

\begin{defn}
A representation of $\FB$ in a $C^*$-algebra $\CC$ is a collection of linear maps $\pi_g : \FB_g \to \CC$ such that
\begin{enumerate}
    \item $\pi_f(\xi_f) \pi_g(\eta_g) = \delta_{\Bs(f), \Bt(g)} \pi_{fg}(\xi_f \eta_g)$,
    \item $\pi_f(\xi_f)^* = \pi_{f^{-1}}(\xi^*_f)$,
\end{enumerate}
for every $f, g \in G(-, -)$, $\xi_f \in \FB_f$, and $\xi_g \in \FB_g$.
\end{defn}

With notations as in \cref{Sub:_fock}, $C_c(\FB)$ is a $*$-algebra with multiplication and involution defined as follows:
\begin{align*}
(\eta \zeta)_h := \sum_{h = fg} \eta_f \zeta_g, \quad (\eta^*)_f := (\eta_{f^{-1}})^*\quad  \forall \eta, \zeta \in C_c(\FB), f, h \in G(-, -).
\end{align*}
We can identify $\xi_f \in \FB_{f}$ as an element in $C_c(\FB)$.
 It is clear there is a one-to-one correspondence between the set of representations of $\FB$ in C$^*$-algebras and the set of  $*$-homomorphisms from $C_c(\FB)$ into C$^*$-algebras. Let $\pi : C_c(\FB) \to \CC$ be a $*$-homomorphism. Note that
\begin{align*}
    \|\pi(\xi_f)\|^2 = \|\pi ((\xi_f)^* \xi_f)\|  \le \| \xi_f^* \xi_f\| = \| \xi_f\|^2.
\end{align*}
 Thus $\|\pi(\xi)\| \leq \sum_{f \in G(-, -)} \| \xi_f \|$ for every  $\xi = (\xi_f) \in C_c(\FB)$.

The \textit{full cross sectional $C^*$-algebra} $C^*(\FB)$ of $\FB$ is the completion of $C_c(\FB)$ with respect to the norm
\begin{align*}
\| \xi \|_u := \sup \{ \| \pi(\xi)\|: \pi \ \ \textit{is a $*$-representation of $C_c(\FB)$}\},
\end{align*}
where $\xi \in C_c(\FB)$. The full cross sectional $C^*$-algebra $C^*(\FB)$ satisfies the following universal property: every $*$-homomorphism $\pi: C_c(\FB) \to \CC$ is factor through $C^*(\FB)$.

If the quasi-product system $\FX$ in  \cref{Sub:_fock} is a Fell bundle $\FB$, then it is not hard to check that $L_{\eta_g}^* = L_{\eta_g^*}$, $\|L_{\eta_g}\| = \|\eta_g\|$ for each $\eta_g \in \FB_g$, $g \in G(-, -)$.
The representation $\{\eta_g \mapsto L_{\eta_g}: \FB_g \to \CL(C_0(\FB))\}_{g \in P(-,-)}$  is called the \textit{regular representation} of $\FB$.
The \textit{reduced cross sectional algebra} $C^*_r(\FB)$ of $\FB$ is the completion of the image of the left regular representation of $\FB$.
\begin{prop}\label{lem:_the_coaction_on_reduced}
    There exists a normal coaction $\delta_{C^*_r(\FB)}$ of $G$ on $C^*_r(\FB)$ such that
    \begin{align*}
        \delta_{C^*_r(\FB)}(L_{\eta_g}) = L_{\eta_g} \otimes U_g, \quad \forall \eta_g \in \FB_g, \quad g \in G(-, -).
    \end{align*}
\end{prop}

\begin{proof}
  By  \cref{prop:_Generalized Fell absorption principle}, we see that there is a $*$-homomorphism $\delta^r_{C^*_r(\FB)} : C^*_r(\FB) \to C^*_r(\FB) \otimes C^*_r(G)$ such that $\delta^r_{C^*_r(\FB)}(L_{\eta_g}) = L_{\eta_g} \otimes \lambda_g$. It is not hard to see that $\delta^r_{C^*_r(\FB)}$ is injective. Note that $\delta^r_{C^*_r(\FB)}$ satisfies the conditions in \cref{pro:_delta_tran_norm}. Thus there exists a normal coaction $\delta_{C^*_r(\FB)}$ of $G$ on $C^*_r(\FB)$ such that $\delta_{C^*_r(\FB)}(L_{\eta_g}) = L_{\eta_g} \otimes U_g$.

\end{proof}

Recall that a \textit{$G$-graded C$^*$-algebra} is a C$^*$-algebra $\CB$ with a $G$-grading, i.e., a family of closed linear independent subspaces $\{\CB_g\}_{g \in G(-, -)}$ 
such that
\begin{enumerate}
\item $\CB_f \CB_g \subseteq \CB_{fg}$, for $\Bs(f) = \Bt(g)$.  Otherwise, $\CB_f \CB_g = 0$;
\item $\CB^*_g = \CB_{g^{-1}}$;
\item $\sum_{f \in G(-, -)}\CB_f$ is dense in $\CB$.
\end{enumerate}

 The  $C^*$-algebras $C^*_r(\FB)$ and $C^*(\FB)$ are examples of $G$-graded C$^*$-algebras. Conversely, by mimicking the proof theorem 3.3 in  \cite{Exel1997}, we have the following result.

\begin{lem}\label{thm:_couniver_groupid_cond}
    Let $\{\CB_g\}_{g \in G(-, -)}$ be a family of closed linear subspaces (not necessarily linear independent) of a C$^*$-algebra $\CB$. Assume that
\begin{itemize}
    \item $\{\CB_g\}$ satisfies the conditions (1)-(3) above,
    \item there is a bounded linear map $F : \CB \to \overline{\spn}\{ b:  b \in \CB_{1_x}, x \in \ob(G) \}$ such that $F|_{\overline{\spn}\{b:  b \in  \CB_{1_x},  x \in \ob(G) \}} = \id$  and $F|_{\CB_g} = 0$  for $g \ne 1_x$, $x \in \ob(G)$.
\end{itemize}
Then we have
\begin{enumerate}
    \item $F$ is a conditional expectation;

    \item $\{\CB_g\}$ is a family of linear independent subspaces, i.e., $\{\CB_g\}$ is a $G$-grading;

    \item $\FB := \{ \CB_g\}_{g\in G(-, -)}$ is a Fell bundle, and there exists a surjective $*$-homomorphism $\lambda_\CB : \CB \to C^*_r(\FB)$ such that $\lambda_{\CB}(\eta_g) = L_{\eta_g}$ for $\eta_g \in \CB_g$, $g \in G(-, -)$.
\end{enumerate}
\end{lem}

\begin{defn}\label{defn:topologically_G_C_alg}
    A topologically $G$-graded $C^*$-algebra is a $G$-graded $C^*$-algebra $\CB$ equipped with a conditional expectation $F$ from $\CB$ to $\overline{\spn}\{b: b \in  \CB_{1_x}, x \in \ob(G) \}$ such that $F|_{\CB_g} = 0$ for $g \neq 1_x$.
\end{defn}

\begin{example}\label{expl:top_G_graded_alg}
    Let $\FB=\{\FB_f\}_{f \in G(-, -)}$ be a Fell bundle over $G$.
    \begin{enumerate}
        \item It is clear that $\FB_f$ can be regarded as a right Hilbert $C_0(\{\FB_{x}\})$-module. Then the canonical inclusion $J_f: \eta_f \mapsto \eta_f: \FB_f \to C_0(\FB)$ is adjointable with $J_f^*(\xi) = \xi_f$, $\xi \in C_0(\FB)$. Let $\zeta \in C_c(\FB)$. We have $J_f^* L_{\zeta} J_{1_{\Bs(f)}} (a)= \zeta_f a$ for every $a \in \FB_{\Bs(f)}$. Therefore, for every $\eta \in C^*_r(\FB)$, there exists a unique $\widehat{\eta}_f \in \FB_f$ such that
        \begin{align}\label{equ:_notation}
        J^*_f \eta J_{1_{\Bs(f)}} (a)= \widehat{\eta}_f a.
        \end{align}
        Moreover, $\|\widehat{\eta}_f\| \leq \|\eta\|$ since $\|J_f\| =1$. By Proposition 3.6 and Proposition 3.10 in \cite{Kumjian1998}, the map defined by
 \begin{align}\label{equ:_conditional_expectation}
     E(\eta) := \sum_{x \in \ob(G)} L_{\widehat{\eta}_{1_x}},
  \end{align}
is a faithful conditional expectation from $C^*_r(\FB)$ onto $C_0(\{\FB_{x}\})$. In particular, the reduced cross sectional algebra $C^*_r(\FB)$ of $\FB$ is a topologically $G$-graded $C^*$-algebra.

    \item Since $\|\xi_f\|_u = \| L_{\xi_f}\| = \| \xi_f\|$ for $\xi_f \in \FB_f$, we can also regard $C_0(\{\FB_x\})$ as a subalgebra of $C^*(\FB)$.
By the universal property of $C^*(\FB)$, there is a canonical $*$-homomorphism $\Lambda_* :  C^*(\FB) \to C^*_r(\FB)$. Let $E_{C^*(\FB)}(-) := E\circ  \Lambda_*(-)$, where $E: C^*_r(\FB) \to C_0(\{\FB_x\})$ is the conditional expectation defined by \cref{equ:_conditional_expectation}. It is obvious that $E_{C^*(\FB)}(-)$ is a conditional expectation from $C^*(\FB)$ onto $C_0(\{\FB_x\})$. Then the full cross sectional $C^*$-algebra $C^*(\FB)$ of $\FB$ is a topologically $G$-graded $C^*$-algebra.


    \end{enumerate}
\end{example}


\begin{cor}\label{cor:_the_existence_of_the_contractive}
Let $\CB$ be a topologically $G$-graded $C^*$-algebra. For each $g \in G(-, -)$, there exists a contractive linear map $F_g : \CB \to \CB_g$ such that $F_g(\xi) = \xi_g$ for finite sums $\xi= \sum_{f \in G(-, -)} \xi_f$, $\xi_f \in \CB_f$.
\end{cor}

\begin{proof}
    By \cref{thm:_couniver_groupid_cond}, $\FB:=\{\CB_g\}$ is a Fell bundle and there exists a $*$-homomorphism $\lambda_\CB: \CB \to C^*_r(\FB)$. Then the contractive linear map is given by $F_g(\xi) := \widehat{\lambda_\CB(\xi)}_g$ defined by \cref{equ:_notation}.
\end{proof}

\begin{cor}\label{pro:_thre_mini_alge_ass_fell}
    Let $\CB$ be a topologically $G$-graded $C^*$-algebra such that the conditional expectation $F: \CB \to \overline{\spn}\{b: b \in  \CB_{1_x}, x \in \ob(G) \}$ is faithful. Then the $*$-homomorphism $\lambda_\CB: \CB \to C^*_r(\FB)$ described in \cref{thm:_couniver_groupid_cond} is injective.
\end{cor}

\begin{proof}
    Let $E: C^*_r(\FB) \to C_0(\{\FB_{x}\})$ be faithful conditional expectation defined by \mref{equ:_conditional_expectation}. Note that $F(b) = E(\lambda_B(b))$ for every $b \in \CB$. Since $F$ is faithful, we have $\lambda_\CB$ is injective.
\end{proof}


\begin{cor}\label{prop:_grading_coaction}
    Let $(\FA, \delta)$ be a cosystem such that $\FA$ is a $C^*$-algebra. Then $\FA$ is a topologically $G$-graded C$^*$-algebra with the $G$-grading $\FA^{\delta}_g :=\{ a \in \FA: \delta(a) = a \otimes U_g\}$, $g \in G(-, -)$.
\end{cor}

\begin{proof}
Let $E_{C^*(G)} : C^*(G) \to \overline{\spn}\{ U_{1_x}: x \in \ob(G)\}$ be the conditional expectation such that $E_{C^*(G)}(U_g) = 0$ for every $g \neq 1_x$ (see \cref{expl:top_G_graded_alg}). Since $\sum_{g \in G(-, -)} \FA^\delta_g $ is dense in $\FA$, we have $(\id \otimes E_{C^*(G)}) \circ \delta (\FA) \subset \delta (\FA)$. Therefore,
\begin{align*}
    E_\FA(-) : = \delta^{-1} \circ (\id \otimes E_{C^*(G)}) \circ \delta (-),
\end{align*}
    is a bounded linear map such that $E_\FA$ is identity on $\overline{\spn}\{ \FA^{\delta}_{1_x}: x \in \ob(G) \}$ and vanishes on each $\FA^\delta_g$, $g \ne 1_x$. By \cref{thm:_couniver_groupid_cond}, $\FA$ is a topologically $G$-graded C$^*$-algebra with the $G$-grading $\{\FA^{\delta}_{g}\}_{g \in G(-,-)}$.
\end{proof}

\section{$C^*$-envelopes of cosystems as co-universal $C^*$-algebras}

The notion of co-universal $C^*$-algebras for  injective, gauge-compatible,  Nica covariant Toeplitz representations of compactly aligned product systems over quasi-lattice ordered groups was introduced by Carlsen, Larsen, Sims and Vittadello in \cite{MR2837016}. The existence of the co-universal algebra for compactly aligned product system over a group-embeddable right LCM semigroup was proved in \cite{DKKLL21}. In this section, we will consider the co-universal algebra of compactly aligned product system over a finite aligned subcategory of a groupoid.

Throughout the rest of this section, we use $\FX = (\{\CA_{x}\}, \{\FX_{f}\})$ to denote a compactly aligned product system over a finite aligned subcategory $P$ of a groupoid $G$, and $(\CB, \delta: \CB \to \CB \otimes C^*(G))$ to denote a cosystem such that $\CB$ is a C$^*$-algebra.  By \cref{rem:_algebra}, $\delta$ is an injective $*$-homomorphism.
\begin{defn}
	A gauge-compatible Toeplitz representation $\psi$ of $\FX$ into $(\CB, \delta)$ is a Toeplitz representation $\psi: \FX \to \CB$ such that
	\begin{align*}
		\delta(\psi_f(\xi)) = \psi_f(\xi) \otimes U_f, \quad  \forall \xi \in \FX_f, f \in P(-, -).
	\end{align*}
\end{defn}

\begin{remark}
	Let $\psi$ be a gauge-compatible Toeplitz representation of $\FX$ into $(\CB, \delta)$. Then the restriction of $\delta$ on $C^*(\psi)$ is a coaction of $G$ on $C^*(\psi)$.
\end{remark}

An injective, gauge-compatible, Nica covariant Toeplitz representation $j^r$ of $\FX$ is said to be co-universal if for every injective, gauge-compatible, Nica covariant Toeplitz representation $\psi$ of $\FX$, there is a (surjective) $*$-homomorphism $\phi: C^*(\psi) \to C^*(j^r)$ such that the following diagram
\begin{align*}
\xymatrix @R=0.2in {
    & \FX \ar[ld]_{\psi} \ar[rd]^-{j^r} &  \\
   C^*(\psi) \ar[rr]_-{\phi}  & & C^*(j^r)
   }
\end{align*}
commutes. Note that $\phi$ is automatically equivariant. The C$^*$-algebra $C^*(j^r)$ is called the co-universal C$^*$-algebra for injective, gauge-compatible, Nica covariant Toeplitz representations of $\FX$.

\begin{lem}\label{pro:_the_coaction_of_TRX}
There exists a normal coaction $\overline{\delta}$ of $G$ on ${\TTT}^r_\FX$ such that
\begin{align}\label{equ:_fock}
	(\TTT^r_\FX)^{\overline{\delta}}_g = \overline{\spn}\{& L_{\xi_{f_1}} L_{\eta_{k_1}}^* \cdots L_{\xi_{f_n}} L_{\eta_{k_n}}^*:  \xi_{f_i} \in \FX_{f_i}, \eta_{k_i} \in \FX_{k_i},f_1 k_1^{-1}  \cdots f_n k^{-1}_n = g, f_i, \nonumber \\
	& k_i \in P(-, -),   1 \leq i \leq n \in \Nb^{+}\},
\end{align}
for every $g \in G(-, -)$.
Moreover, ${\TTT}^r_\FX$ is a topologically $G$-graded $C^*$-algebra and
\begin{align*}
C^*_r(\{(\TTT^r_\FX)^{\overline{\delta}}_g \}_{g \in G(-, -)}) \cong \TTT^r_\FX.
\end{align*}
\end{lem}

\begin{proof}

With notations as in \cref{Sub:_fock}.  By \cref{prop:_Generalized Fell absorption principle}, there is a $*$-homomorphism  $\overline{\delta}_r: {\TTT}^r_\FX \to {\TTT}^r_\FX \otimes C^*_r(G)$ such that $\overline{\delta}_r(L_{\eta_k}) = L_{\eta_k} \otimes \lambda_k$, for $k \in P(-, -)$. It is not hard to check that $\overline{\delta}_r$ is injective.
    By \cref{pro:_delta_tran_norm}, there exists a normal coaction $\overline{\delta}$ of $G$ on ${\TTT}^r_\FX$ such that $\overline{\delta}(L_{\eta_k}  ) = L_{\eta_k} \otimes U_{k}$.    In particular,
    \begin{align*}
    	\overline{\spn}\{ L_{\xi_{f_1}} L_{\eta_{k_1}}^* \cdots L_{\xi_{f_n}} L_{\eta_{k_n}}^*:&\ \  \xi_{f_i} \in \FX_{f_i}, \eta_{k_i} \in \FX_{k_i},f_1 k_1^{-1}  \cdots f_n k^{-1}_n = g, f_i, \\
    	& k_i \in P(-, -),   1 \leq i \leq  n \in \Nb^{+}\} \subseteq (\TTT^r_\FX)^{\overline{\delta}}_g,
    \end{align*}
    since $\overline{\delta}$ is a $*$-homomorphism.

Let $b \in ({\TTT}^r_\FX)^{\overline{\delta}}_g$ and $\epsilon >0$. Then there exists
\begin{align*}
	a \in {\spn}\{L_{\xi_{f_1}} L_{\eta_{k_1}}^* \cdots L_{\xi_{f_n}} L^*_{\eta_{k_n}} :&\ \ \xi_{f_i} \in \FX_{f_i},  \eta_{k_i} \in \FX_{k_i}, f_i, k_i \in P(-,-),\\
	&1 \leq i \leq n \in \Nb^{+}\}
\end{align*}
such that $\|b - a\| \leq \epsilon$. Recall that $C^*(G)$ is a topologically $G$-graded C$^*$-algebra (see \cref{expl:top_G_graded_alg}). Thus there is a contractive linear map $P_{g}$ from $C^*(G)$ onto $\Cb U_{g}$ such that $P_g(\sum_{h \in G(-, -)} a_h U_h) = a_g U_g$, for every $g \in G(-, -)$ (see \cref{cor:_the_existence_of_the_contractive}). Note that
\begin{align*}
	\|(\id  \otimes P_{g})(b \otimes U_g - \overline{\delta}(a))\| \leq  \|b \otimes U_g - \overline{\delta}(a)\| = \|b - a\| \leq \epsilon.
\end{align*}
Therefore, we may assume that
\begin{align*}
	a \in \overline{\spn}\{ L_{\xi_{f_1}} L_{\eta_{k_1}}^* \cdots L_{\xi_{f_n}} L_{\eta_{k_n}}^*:& \ \ \xi_{f_i} \in \FX_{f_i}, \eta_{k_i} \in \FX_{k_i},f_1 k_1^{-1}  \cdots f_n k^{-1}_n = g, \\
	&f_i, k_i \in P(-, -),   1 \leq i \leq n \in \Nb^{+}\}.
\end{align*}
 Thus \cref{equ:_fock} is hold.

By \cref{prop:_grading_coaction},  ${\TTT}^r_\FX$ is a topologically $G$-graded $C^*$-algebra with  a conditional expectation $E : \TTT^r_\FX \to \overline{\spn}\{  \cup_{x \in \ob(P)} (\TTT_\FX^r)^{\bar{\delta}}_{1_x}\}$. For every $f \in P(-, -)$, let $Q_f$ be the projection from $C_0(\FX)$ onto $\FX_f$. It is not hard to see that
    \begin{align}\label{alignequ:full_reduced_algebra_of_NX_1}
     E (b) = \sum_{f \in P(-, -)} Q_f b Q_{f}, \quad \forall  b \in {\TTT}^r_\FX.
    \end{align}
 Thus $E$ is a faithful conditional expectation. Then by  \cref{pro:_thre_mini_alge_ass_fell}, $\TTT^r_\FX$ is $*$-isomorphic to $C^*_r(\{(\TTT^r_\FX)^{\overline{\delta}}_g \}_{g \in G(-, -)})$.
\end{proof}

Let $\psi: \FX \to \CB(\FH)$ be a Nica covariant Toeplitz representation of $\FX$. For every finite subset $\F \subseteq P(-, -)$, let
\begin{align*}
	\psi(\F):=
	\begin{cases}
		\Pi_{f \in \F} \psi^{(f)}(I), & \F \neq \emptyset;\\
		I, & \F = \emptyset.
	\end{cases}
\end{align*}
Note that
\begin{align}\label{eq:_the_linea_sp}
	B_{\G, \psi} := \overline{\spn}\{&\psi^{(f_1)}(S_{f_1}) \cdots \psi^{(f_n)}(S_{f_n})\psi(\F): \\
	&S_{f_i} \in \CK(\FX_{f_i}), f_i \in \F,   \F \mbox{ is a finite subset of $\G$}, n \geq 1\}\nonumber
\end{align}
is a C$^*$-subalgebra of $[C^*(\psi)]$ (see \cref{gene:algebra}), for every subset $\G$ of $P(-, -)$. It is clear that
\begin{align}\label{equ:_the_core_euqavar}
	B_{P, \psi} = \overline{ \cup \{B_{\F,\psi} : \mbox{$\F$ is a finite subset of $P(-, -)$ }\}}.
\end{align}
By \cref{lem:_incluP}, \cref{lem:the_op} and \cref{remarkpro}, we have
\begin{align*}
	B_{P, \psi} = \overline{\spn}\{&\psi ^{(f)}(S_f) \psi(\F): \\
	&S_{f} \in \CK(\FX_{f}), f \in \F, \text{ $\F$ is a finite subset of $P(-, -)$} \}.
\end{align*}

Recall that $\psi$ can be decomposed as a direct sum of cyclic Nica covariant Toeplitz representations. It is easy to see that there is a $*$-homomorphism $\Phi_{\psi} : [\NNN\TTT_\FX] \to [C^*(\psi)]$ such that
\begin{align*}
	\Phi_{\psi}(i_{\FX,f}(\xi_f) i_{\FX, g}(\eta_g)^* i_\FX(\F)) = \psi_f(\xi_f) \psi_g(\eta_g)^* \psi(\F),
\end{align*}
for all $\xi_f \in \FX_f$, $\eta_g \in \FX_g$, $f, g \in P(-, -)$, $\F$ is a finite subset of $P(-, -)$.


With the notation as in \cref{lem:Fockr}, we use $\widetilde{\TTT}^r_\FX$ to denote the $C^*$-algebra generated by the image of $\widetilde{L}$. Note that $\TTT_X^r \simeq \widetilde{\TTT}_X^r$.

\begin{lem}\label{lem:_the_fell_bundle_isomorphism}
	With notations as above. The restriction of $\Phi_{\widetilde{L}}: [\NNN\TTT_\FX] \to [\widetilde{ \TTT}_\FX^r]$ on $B_{P, i_\FX}$ is injective.
\end{lem}

\begin{proof}
	Let $\F$ be a finite subset of $P(-, -)$. By \cref{equ:_the_core_euqavar}, we only need to show that $\Phi_{\widetilde{L}}$ is injective on
	\begin{align*}
		B_{\F, i_\FX} &:= \overline{\spn}\{i_\FX^{(f_1)}(S_{f_1}) \cdots i_\FX^{(f_n)}(S_{f_n})i_\FX(\G): \\
		&S_{f_i} \in \CK(\FX_{f_i}), f_i \in \G,   \G \mbox{ is a finite subset of $\F$}, n \geq 1\}.
	\end{align*}
	In the following, we use $|\F|$ to denote the cardinality of $\F$.
	
	Note that
	\begin{align*}
		\vee_{k \in \F} i_{\FX}^{(k)}(I) = \sum_{\emptyset \neq \F' \subseteq \F} E_{\F'},
	\end{align*}
	where
	\begin{align*}
		E_{\F'} := \sum_{\F' \subseteq \F'' \subseteq \F} (-1)^{|\F''|- | \F'|} i_\FX(\F'') = i_\FX(\F') \left (I - \vee_{f \in \F \setminus \F'} i_{\FX}^{(f)}(I) \right).
	\end{align*}
	Similarly, let
	\begin{align*}
		\widetilde{E}_{\F'} := \sum_{\F' \subseteq \F'' \subseteq \F} (-1)^{|\F''|- | \F'|} \widetilde{L} (\F'') = \widetilde{L}(\F') \left (I - \vee_{f \in \F \setminus \F'} \widetilde{L}^{(f)}(I) \right).
	\end{align*}
	It is not hard to check that $\widetilde{E}_{\F'} \neq 0$ if $\left (\cap_{f \in \F'} f P \right) \setminus \cup_{g \in \F \setminus \F'} gP \neq \emptyset$, and $E_{\F'} = 0$ if $\left (\cap_{f \in \F'} f P \right) \setminus \cup_{g \in \F \setminus \F'} gP = \emptyset$.
	Thus $\widetilde{E}_{\F'} \neq 0$ if and only if  $E_{\F'} \neq 0$.
	
	Let $\{b_j = \sum_{\emptyset \neq \G \subseteq \F} b_{ \G, j} i_\FX(\G)\}_j$ be a sequence in $B_{\F, i_\FX}$, where $b_{ \G, j}$ is in
	\begin{align*}
		\spn\{i_\FX^{(f_1)}(S_{f_1}) \cdots  i_\FX^{(f_n)}(S_{f_n}): S_{f_i} \in \CK(\FX_{f_i}), f_i \in \G, n \geq 1 \}.
	\end{align*}
	Assume that $b_j$ converges to $b$ such that $\Phi_{\widetilde{L}}(b) = 0$.
	
	We claim that $b_j E_{\F'}$ converges to $0$ when
	$\left (\cap_{f \in \F'} f P \right) \setminus \cup_{g \in \F \setminus \F'} gP \neq \emptyset$.
	 Indeed, let $\omega \in \vee \{f: f \in \F'\}$ such that $\omega \notin \cup_{g \in \F \setminus \F'} gP$. We only need to show that $b_j i_\FX^{(\omega)}(I)\left (I - \vee_{f \in \F \setminus \F'} i_{\FX}^{(f)}(I) \right)$ converges to $0$.
	
	Note that
	\begin{align*}
		b_j &= \sum_{\emptyset \neq \F' \subseteq \F} \left(\sum_{\emptyset \neq \G \subseteq \F'}  b_{\G, j} \right) E_{\F'}.
	\end{align*}
	By \cref{lem_the_reg_lift}, there exists $c_{\G, j} \in \CL(\FX_\omega)$ such that
	\begin{align*}
		b_j i_\FX^{(\omega)}(I) \left (I - \vee_{f \in \F \setminus \F'} i_{\FX}^{(f)}(I) \right) &= i_{\FX}^{(\omega)}\left(\sum_{\emptyset \neq \G \subseteq \F'}  c_{\G, j} \right) \left (I - \vee_{f \in \F \setminus \F'} i_{\FX}^{(f)}(I) \right),\\
		\Phi_{\widetilde{L}}\left (b_j i_\FX^{(\omega)}(I)  \left (I - \vee_{f \in \F \setminus \F'} i_{\FX}^{(f)}(I) \right)\right)&= \widetilde{L}^{(\omega)}\left(\sum_{\emptyset \neq \G \subseteq \F'}  c_{\G, j} \right)  \left (I - \vee_{f \in \F \setminus \F'} \widetilde{L}^{(f)}(I) \right).
	\end{align*}
	Since $\Phi_{\widetilde{L}}\left (b_j i_\FX^{(\omega)}(I)  \left (I - \vee_{f \in \F \setminus \F'} i_{\FX}^{(f)}(I) \right)\right)$ converges to $0$ and
	 \begin{align*}
	 \|\widetilde{L}^{(\omega)}\left(\sum_{\emptyset \neq \G \subseteq \F'}  c_{\G, j} \right)  \left (I - \vee_{f \in \F \setminus \F'} \widetilde{L}^{(f)}(I) \right)\| = \| \sum_{\emptyset \neq \G \subseteq \F'}  c_{\G, j}\|,
	 \end{align*} we have $\sum_{\emptyset \neq \G \subseteq \F'}  c_{\G, j}$ converges to $0$. Thus $b_j i_\FX^{(\omega)}(I)\left (I - \vee_{f \in \F \setminus \F'} i_{\FX}^{(f)}(I) \right)$ converges to $0$. This completes the proof.
\end{proof}

We are now ready to prove the main result of the paper.

\begin{proof}[\bf{Proof of Theorem 1}]
	Let $\Big(C^*_{e}\big((\TTT^r_\FX)^+ \big), \iota: (\TTT^r_\FX)^+ \to C^*_{e}\big((\TTT^r_\FX)^+\big)\Big)$ be the C$^*$-envelope of $(\TTT^r_\FX)^+$. By the definition of C$^*$-envelope, the completely isometric homomorphism $\iota: (\TTT^r_\FX)^+ \to C^*_{e}\big((\TTT^r_\FX)^+\big)$ can be promoted to a $*$-homomorphism, still denoted by $\iota$, from $\TTT^r_\FX$ to $C^*_{e}\big((\TTT^r_\FX)^+\big)$.
	
	Let $\overline{\delta}: \TTT^r_\FX \to \TTT^r_\FX \otimes C^*(G)$ be the normal coaction of $G$ on $\TTT^r_\FX$ defined in \cref{pro:_the_coaction_of_TRX}. It is clear that the restriction of $\overline{\delta}$ on $(\TTT^r_\FX)^+$ is a normal coaction of $G$ on $(\TTT^r_\FX)^+$. By \cref{lem:_envelope_for_the_cosystem}  and  the proof of \cref{pro:_the_coaction_of_TRX}, the C$^*$-envelope of the cosystem $\big((\TTT^r_\FX)^{+}, \overline{\delta}|_{(\TTT^r_\FX)^+}\big)$ is
	\begin{align*}
		\Big(C^*\big(\{ j^r_{f}(\xi_f) : \xi_f \in \FX_f, f \in P(-, -)\}\big),\id \otimes \delta_G, (\iota \otimes \id_{C^*(G)}) \overline{\delta}|_{(\TTT^r_\FX)^+}\Big),
	\end{align*}
	where $j^r_{f}: \FX_f \to C^*_{e}\big((\TTT^r_\FX)^+\big) \otimes C^*(G)$ is the linear map defined by $j^r_{f}(\xi_f) := \iota (L_{\xi_f} )\otimes U_f$. Recall that $L$ is injective and Nica covariant, we have $\{j^r_{f}\}_{f \in P(-,-)}$ is an injective, gauge-compatible, Nica covariant Toeplitz representation of $\FX$ by the proof of \cref{Nica_cond}. We show next that $\{j^r_f\}_{f \in P(-,-)}$ is co-universal.
	
	First note that, by \cref{lem:the_op}, \cref{lem:_incluP} and \cref{remarkpro}, we have
	\begin{align*}
		\widetilde{B}_{P, i_\FX} &:= \overline{\spn}\{ i_{\FX, f_1}(\xi_{f_1}) i_{\FX, k_1}(\eta_{k_1})^*  \cdots  i_{\FX, f_n}(\xi_{f_n}) i_{\FX, k_n}(\eta_{k_n})^* :   \xi_{f_i} \in \FX_{f_i}, \eta_{k_i} \in \FX_{k_i},\\
		&f_1 k_1^{-1}  \cdots f_n k^{-1}_n = 1_{\Bt(f_1)}, f_i,  k_i \in P(-, -),   1 \leq i \leq n \in \Nb^{+}\} \subseteq B_{P, i_\FX}.
	\end{align*}
	Then the canonical $*$-homomorphism $\Phi$ from $\NNN\TTT_\FX$ onto $\TTT^r_\FX$ is injective on $\widetilde{B}_{P, i_\FX}$ by \cref{lem:_the_fell_bundle_isomorphism}.

Let $\psi$ be an injective, gauge-compatible,  Nica covariant Toeplitz representation of $\FX$. Then $C^*(\psi)$ is a topologically $G$-graded $C^*$-algebra by \cref{prop:_grading_coaction}. We use $\FB$ to denote the $G$-grading of $C^*(\psi)$. Note that $\FB$ is a Fell bundle. Let
\begin{align*}
	\phi := (\NNN\TTT_\FX \to C^*(\psi) \to C^*_r(\FB)),
\end{align*}
where the first arrow is the canonical $*$-homomorphism from $\NNN\TTT_\FX$ to $C^*(\psi)$, and the second arrow is the surjective $*$-homomorphism from $C^*(\psi)$ to the reduced cross sectional $C^*$-algebra $C^*_r(\FB)$ of $\FB$ given by \cref{thm:_couniver_groupid_cond}.

 Let $\widehat{E}: \NNN\TTT_\FX \to \widetilde{B}_{P, i_\FX}$ be the expectation defined by
\begin{align*}
	\widehat{E}(-) :=  \Phi|_{\widetilde{B}_{P, i_\FX}}^{-1} \circ \overline{E} \circ \Phi(-),
\end{align*}
where $\overline{E}$ is the expectation of $\TTT_\FX^r$ defined as in the proof of \cref{pro:_the_coaction_of_TRX}.
For every $b \in \ker(\Phi)$, we have
\begin{align*}
	E' \phi(b^*b) = \phi \widehat{E}(b^*b) = 0,
\end{align*}
where $E' : C^*_r(\FB) \to  C_0(\{\FB_x\})$ is the faithful conditional expectation defined by \cref{equ:_conditional_expectation}. Thus, $\phi(\ker(\Phi)) = 0$ and there is a  $*$-homomorphism $\phi'$ such that the following diagram
\begin{align*}
	\xymatrix{
		\NNN\TTT_\FX \ar[r]^{\Phi} \ar[rd]_{\phi}  & \TTT^r_\FX \ar[d]^-{\phi'} \\
		& C^*_r(\FB)
	}
\end{align*}
commutes.

By \cref{lem:_the_coaction_on_reduced}, there exists a normal coaction $\delta_{C^*_r(\FB)}: C^*_r(\FB) \to C^*_r(\FB) \otimes C^*(G)$.
Since $\phi': \TTT^r_\FX \to C^*_r(\FB)$ is a $*$-homomorphism, $\phi'$ is completely contractive. Let $\{E_{i,j}\}_{i,j \in n}$ be the caonical matrix units of $M_n(\Cb)$. Recall that $\phi'_n : (\TTT^r_\FX)^{+} \otimes M_n(\Cb) \subset \CL(C_0(\FX)\otimes \Cb^{(n)})\to C^*_r(\FB) \otimes M_n(\Cb)\subset \CL(C_0(\FB) \otimes \Cb^{(n)})$ is defined by
\begin{align*}
	\phi'_n(\sum_{i,j}(a_{i,j} \otimes E_{ij})) = \sum_{i,j} (\phi'(a_{i, j}) \otimes E_{i, j}),
\end{align*}
where $a_{i,j} \in (\TTT^r_\FX)^{+}$, $1 \leq i, j\leq n \in \Nb^{+}$.
Since $\psi_{1_x}$ is injective for every $x \in \ob(G)$, we can regard $\FX_f$ as a subspace of $\FB_f$ for $f \in P(-, -)$.  It follows that we can regard  $C_0(\FX)$ as a subspace of $C_0(\FB)$. Then we can regard $b \in (\TTT^r_\FX)^{+}$ as the operator $\phi'(b)|_{C_0(\FX)}$. Then we have
\begin{align*}
	\left \|\phi_n' \left (\sum_{i,j} a_{i, j}\otimes E_{i,j} \right ) \right \|  \geq   \left \| \sum_{i,j}(a_{i,j} \otimes E_{i,j}) \right \|.
\end{align*}
Thus $ \phi'|_{(\TTT^r_\FX)^+}$ is completely isometric. It is not hard to check that $ \phi'$ is equivariant. Then $$(C^*_r(\FB),  \delta_{C^*_r(\FB)}, \phi'|_{(\TTT^r_\FX)^+})$$ is a $C^*$-cover of the cosystem $((\TTT^r_\FX)^+,  \overline{\delta})$.

Therefore, there is a $*$-homomorphism
\begin{align*}
	\phi'' : C^*(\psi) \to C^*_r(\FB) \to C^*(\{ j^r_{f}(\xi_f): \xi_f \in \FX_f, f \in P(-, -)\})
\end{align*}
such that $\phi'' \psi(\xi_f) = j^r_{f}(\xi_f)$, for every $\xi_f \in \FX_f$ and $f \in P(-, -)$, and $\{j^r_{f}\}_{f \in P(-,-)}$ is co-universal. Finally, by the discussion after \cref{lem:_envelope_for_the_cosystem}, and the proof of  \cref{pro:_the_coaction_of_TRX}, we have the coaction of $G$ on the co-universal C$^*$-algebra is normal.
\end{proof}

\bibliographystyle{amsplain}

\end{document}